\documentclass{article}
\usepackage[english]{babel}
\usepackage{authblk}
\usepackage[letterpaper,top=2cm,bottom=2cm,left=3cm,right=3cm,marginparwidth=1.75cm]{geometry}
\usepackage{amsmath,amsfonts,amsthm,authblk}
\usepackage{graphicx}
\usepackage{xcolor}
\usepackage[colorlinks=true, allcolors=blue]{hyperref}

\newcommand{\p}{\partial}
\newcommand{\R}{\mathbb{R}}
\newcommand{\NN}{\mathbb{N}}

\newcommand{\CM}{\mathcal{M}}
\newcommand{\FM}{\mathrm{FM}}

\newcommand{\BL}{\mathrm{BL}}
\newcommand{\TV}{\mathrm{TV}}

\newcommand{\ind}{{1\!\!1}}
\newcommand{\eps}{\varepsilon}

\newtheorem{theorem}{Theorem}[section]
\newtheorem{lemma}{Lemma}[section]
\newtheorem{prop}{Proposition}[section]
\newtheorem{clry}{Corollary}[section]
\newtheorem{definition}{Definition}[section]
\newtheorem{remark}{Remark}[section]

\newtheorem*{remark*}{Remark}
\usepackage{hyperref}
\usepackage{xcolor}

\title{Invariance properties of the solution operator for measure-valued semilinear transport equations}
\author{Sander C. Hille$^1$, Rainey Lyons$^2$, Adrian Muntean$^3$}
\date{
$^1$ Mathematical Institute, Leiden University, The Netherlands\\%
$^2$ Department of Applied Mathematics, University of Colorado Boulder, USA\\%
$^3$Department of Mathematics and Computer Science, Karlstad University, Sweden\\[2ex]%
\today}

\begin{document}

\maketitle

\begin{abstract}
     We provide conditions under which we prove for measure-valued transport equations with non-linear reaction term in the space of finite signed Radon measures, that positivity is preserved, as well as absolute continuity with respect to Lebesgue measure, if the initial condition has that property. Moreover, if the initial condition has $L^p$ regular density, then the solution has the same property. 
\end{abstract}


\section{Introduction}

We consider the semilinear measure-valued transport equation with reaction term
\begin{equation}\label{Eq:Problem_P}\tag{$P$}
\left \lbrace
\begin{split}
    \p_t \mu_t + \text{div}_x(v_ t(x) \mu_t) &= f_t(\mu_t)  &\quad (t,x) \in (0,T)\times \Omega =: \Omega_T,\\
    \mu_{0} &= \nu\in\CM(\Omega). &
\end{split}\right.
\end{equation}
Here and throughout the manuscript, the space $\Omega$ is taken to be either $\R^d$ or the compact space $\mathbb{T}^d =\R^d/\mathbb{Z}^d$, which is the $d$-dimensional torus or, equivalently, the unit hypercube $[0,1]^d$ with periodic boundary conditions. 
These particular choices of $\Omega$ indicate that we do not include any boundary effects in our current analysis. 
For discussions around this class of transport equations with particular boundary conditions, we refer the reader, for instance, to \cite{evers2014_CRM,Evers_Thesis,evers2015_JDE,evers2016measure}.
The reaction map $(t,\mu)\mapsto f_t(\mu)$ takes values in the vector space $\CM(\Omega)$ of finite signed Radon measures on $\Omega$ and is possibly non-linear and satisfies specific structural conditions given below. 
Well-posedness for equations of this type in the space of Radon measures on $\R^d$ has been studied in \cite{Canizo_ea:2012} for the case of a time-independent velocity field and in \cite{ackleh2020well} for the time-dependent case and we will employ conditions on $v$ and $f$ similar to these works. These include a (local) Lipschitz condition on $f_t$ in the Fortet-Mourier or flat metric on $\CM(\Omega)$.

Our focus in this manuscript is establishing `invariance properties' of the solution operator to the Cauchy problem \eqref{Eq:Problem_P} that maps the initial condition $\nu$ to the unique maximal mild solution $t \mapsto \mu_t$ with initial condition $\mu_0 = \nu$. Let $[0, T^*_\nu)$ be the corresponding maximal interval of existence.
Specifically, starting from the previously mentioned well-posedness results for Problem \eqref{Eq:Problem_P} in the space $\CM(\Omega)$, we provide conditions on the reaction term $f$ such that the solution operator preserves:
\begin{enumerate}
    \item {\bf Positivity}, i.e., a positive initial condition leads to a positive  solution $\mu_t \in \CM(\Omega)^+$ for all $t \geq 0$;
    \item {\bf $L^p$-density}, i.e., an initial condition with $L^p$-density ($1<p \leq\infty$) with respect to the Lebesgue measure $\lambda$ leads to a solution with $\mu_t$ absolutely continuous solution with respect to Lebesgue measure $\lambda$ with $L^p$-density, for $t$ possibly in a proper subinterval $[0,\hat{T}^*_\nu)$ of the maximal interval of existence.
\end{enumerate}
We take the time to point out that these preserved properties are important for various applications of similar reaction-transport equations set in the space of measures, including the study of structured population models \cite{rainey2020,ackleh2021structured,dull2021spaces,jablonski2019measure}, traffic flow \cite{bressan2011optima,bressan2020traffic} 
,  and cell migration \cite{colombi2015differentiated}.

Because of this `measure-centric' line of argumentation, we limit attention to `classical' regular velocity fields $v_t$, which are bounded and globally Lipschitz, allowing us to focus on the technicalities involved in showing that results in the function space (density) setting can be derived from results in the general measure setting. Various auxiliary results that we obtained for this purpose seem to be new (in particular Theorem \ref{Lem:BoundedWeakLimits}) and we believe them to be of separate interest in other contexts too.
It is also interesting to consider how our results hold for transport equations and conservation laws defined by less regular velocity fields, see e.g., \cite{DiPerna-Lions:1989,Bogachev-Mayer-Wolf:1999,Ambrosio-Crippa:2014}, but we leave this for future work.

Study of Problem (P) often occurs as preparation for studying the Cauchy problem for the quasi-linear equation
\begin{equation}\label{eq:quasi-linear problem}
    \partial_t \mu_t + \mathrm{div}_x(v_t[\mu_t]\mu_t) = f_t(\mu_t),\qquad (t,x)\in\Omega_T
\end{equation}
in the space of measures. Various authors have considered well-posedness of Cauchy Problem (P) or \eqref{eq:quasi-linear problem} in the space of Radon measures $\CM(\Omega)$ on $\Omega$ as here \cite{ackleh2020well, Canizo_ea:2012}, or on an interval with boundary interaction \cite{evers2014_CRM, evers2015_JDE, evers2016measure}. For the semilinear problem, the solution concept of focus will be that of {\it mild solutions}. 
These are weak solutions \cite{Evers_MildAreWeak} as commonly considered. 
However, a mild solution is such that $t\mapsto\mu_t$ is continuous when the Radon measures are equipped with the Fortet-Mourier metric. Hence, the solution operator, mapping initial value to the unique associated solution $\mu_\bullet$, defines a dynamical system in the space of Radon measures with continuous motions -- as is usually considered.

We want solutions to possess properties that make the model meaningful to the application of interest, e.g., positivity is required if the measure represents distribution of mass or individuals. 
Often we find it convenient to have densities described relative to a preferred reference measure. Then we like to have that regularity of the density with respect to this reference is preserved. This latter property cannot always be guaranteed in this model, as we will show. Other model equations, as for instance those used in population dynamics, often allow the possibility of aggregation in finite time, like the ones describing chemotaxis. This feature naturally leads to the blow-up in finite time of the classical, weak, or strong concepts of solution  \cite{Winkler}. After this blow-up time,  the solution usually turns into a Dirac -- possibly continuing to move in space (cf. e.g. \cite{Velazquez, Schmeiser}). Hence, the overall framework fits best to a setup in terms of measure solutions. 

The main results in this paper can be summarized as follows -- informally -- because the concepts and notation required for making a precise statement still need introduction.
\begin{itemize}
    \item {\bf Theorem \ref{Lem:BoundedWeakLimits}}: 
    Concerns a result that holds in general Polish spaces $S$. It states, that if one has a sequence of positive finite measures $\mu_n$ that converges to a limit measure $\mu_*$ in the $C_b(S)$-weak topology on measures (or equivalently the Fortet-Mourier or flat metric) such that all $\mu_n$ are absolutely continuous with respect to a single positive Radon measure $\nu$ that is not necessarily finite, with density $f_n$ that is an $L^p$-function for $\nu$ with $1<p\leq\infty$; then, $M=\liminf_{n\to\infty} \|f_n\|_p<\infty$ is sufficient to conclude that $\mu_*$ is also absolutely continuous with respect to $\nu$ and has density $f\in L^p(\nu)$, with $\|f\|_p\leq M$. In particular this holds when $\sup_n\|f_n\|_p\leq M$.
    
    \item {\bf Theorem \ref{thrm:positivity}}: Let the reaction map $f_t(\mu)$ be the sum of a positive production term  $p_t(\mu)$ and a term that is absolutely continuous with respect to $\mu$ and that hence represents degradation (death) or proliferation. That is, 
    \begin{equation}\label{eq:intro shape of f}
        f_t(\mu)=p_t(\mu)  + F_t(\mu)\mu\qquad \mbox{for all}\ t\geq 0,\ \mu\in\CM(\Omega),
    \end{equation} 
    where  $F_t(\mu)$ a bounded Borel-measurable function that can be interpreted as local rate of change. Assume that $F_t(\mu)$ is locally bounded, in the sense that for every $R>0$ and $T>0$, there exists $C'_{R,T}>0$ such that  
    \begin{equation}
        \sup_{0\leq t\leq T}\|F_t(\mu)\|_\infty\leq C'_{R,T} \qquad \mbox{for all}\ \mu\in\CM^+(\Omega)\ \mbox{with}\ \|\mu\|_\TV=\mu(S)\leq R.
    \end{equation}
    Then, the mild solution $\mu_\bullet$ to Problem (P) with initial condition $\mu_{0}=\nu\in\CM^+(\Omega)$ is positive for all its time of existence.
    
    \item {\bf Theorem \ref{thrm:Lp invariance positive}}: Let $f_\bullet$ again be of the form \eqref{eq:intro shape of f} discussed above. Moreover, assume that for some $1<p\leq\infty$, $f_t$ maps $L^p(\lambda)$ into itself for all $t\geq 0$ and is an $L^p$-bounded map, locally uniformly in time, i.e. for any $t>t'$ and $r>0$, there exists $B^f_{t',t}(r)\geq 0$  such that
    \begin{equation}
         \bigl\| f_s(\phi\, d\lambda) \bigr\|_{L^p(\lambda)}\ \leq\ B^f_{t',t}(r)\quad\mbox{for all}\ \phi\in L^p(\lambda),\ \|\phi\|_{L^p(\lambda)} \leq r,\ t'\leq s\leq t.
    \end{equation}
    If $\nu\in\CM^+(\Omega)$ is a (positive) initial condition that is absolutely continuous with respect to Borel-Lebesgue measure $\lambda$ on $\Omega$, with  $L^p$-density, then the maximal mild measure-valued solution $\mu_\bullet$ to Problem (P), defined on $[0,T^*_\nu)$, will have $\mu_t$ absolutely continuous with respect to $
    \lambda$ for all $t$ in a subinterval $[0,\hat{T}^*_\nu)$, with $L^p$-density $\phi_t\in L^p(\lambda)$. The function $t\mapsto \|\phi_t\|_{L^p(\lambda)}$ is bounded on any $[0,T]\subset [0,\hat{T}^*_\nu)$ and $\hat{T}^*_\nu$, which may be strictly smaller than the existence time $T^*_\nu$ of the mild solution in measures, is the time of blow-up of the $L^p$-norm of the $\phi_t$.
    
\end{itemize}
Theorem \ref{Lem:BoundedWeakLimits} is a structural result for general measures in $\CM^+(S)$, while  Theorem \ref{thrm:positivity} and Theorem \ref{thrm:Lp invariance positive} are specific to the class of quasilinear evolution equations in measures where Problem (P) belongs. These results clarify sufficient conditions that ensure the preservation of positivity and of $L^p$ integrability of suitable (mild) solutions. It is worth noting that these properties are generally wanted, but they do not always hold true. We clarified the matter for the case of the scalar equation \eqref{eq:quasi-linear problem}, but our results can be generalized for the case of systems and regardless of the space dimensions, if the coupling in the model components fits our working assumptions. 

In this paper, our goal is to prove in full details  these three main results. We do so in Section \ref{invariance}, after introducing the needed preliminary results in Section \ref{preliminaries} as well as presenting a toolbox in Section \ref{weaks} about working with weak limits of sequences of measures with $L^p$-densities.

\section{Preliminaries}\label{preliminaries}

We begin by describing the notation which will be used throughout the manuscript.
We denote the spatial domain of our problem of interest (P) by $\Omega$, which is either $\R^d$ or the $d$-dimensional torus $\mathbb{T}^d:=\R^d/\mathbb{Z}^d$, equipped with the standard Euclidean-based metrics.
While our focus is on problem \eqref{Eq:Problem_P}, various of the forthcoming results hold for measures over more general Polish spaces, i.e. separable completely metrizable topological spaces. These are of recent interest in the study of structured populations dynamics where the structure variable is represented in by elements in a graph or measure space \cite{dull2024structured}.
To indicate such a setting clearly, we will denote a general Polish space by $S$. 
Any metric, $d$, which metrizes the topology of $S$ is called {\it admissible} for $S$ and we will implicitly assume $S$ is equipped with some admissible metric. 
Occasionally, we shall use $S$ or $S'$ to denote a general measurable space.

\subsection{Function spaces, measures, distances and norms}

We employ common notation for function spaces throughout the manuscript. 
For instance, we let $B_b(S)$ denote the vector space of real-valued bounded Borel measurable functions on $S$ and $C_b(S)$ its subspace of continuous bounded functions. Both are equipped with the supremum norm $\|\cdot\|_\infty$.
The space $\BL(S) \subset C_b(S)$ consists of all bounded Lipschitz functions on $S$ for a (fixed, implicit) admissible metric $d$. 
We will use $|f|_L$ to denote the Lipschitz constant of $f\in\BL(S)$. $\BL(S)$ is a Banach space for the norm $\|f\|_\FM := \max\bigl(\|f\|_\infty, |f|_L\bigr)$. One has $\|fg\|_\FM\leq 2\|f\|_\FM\|g\|_\FM$ for $f,g\in\BL(S)$.

We denote by $\mathcal{B}(X)$ the Borel $\sigma$-algebra on the topological space $X$.
$\CM(S)$ represents the vector space of finite signed Borel measures on $S$ and $\CM^+(S)$ the convex cone of positive measures in $\CM(S)$. 
Since $S$ is Polish, every finite Borel measure is a finite Radon measure (compare \cite{BogachevII:2007}). 
Various definitions of Radon measures appear throughout the literature, all leading to the same concept for finite measures (cf. e.g. \cite[Section 7.14(xviii)]{BogachevII:2007} or \cite{Folland:1999}).
We follow \cite[Section 7.1, p.212]{Folland:1999} by defining a positive (possibly infinite) Radon measure to be a Borel measure $\mu$ that is (a) finite on compact sets, (b) outer (open) regular on all Borel sets, i.e.
\begin{equation}\label{eq:outer open regular}
    \mu(E) = \inf\bigl\{ \mu(U)\colon U\subset S\ \mbox{open},\ E\subset U\bigr\},\qquad E\subset S\ \mbox{Borel}
\end{equation}
and (c) inner (compact) regular on all open sets:
\[
    \mu(U) = \sup\{\mu(K)\colon K\subset S\ \mbox{compact}, K\subset U\bigr\},\qquad U\subset S\ \mbox{open}.
\]
In particular, the Borel-Lebesgue measure $\lambda$ on $\R^d$ is a ($\sigma$-finite) Radon measure in this sense. We shall use the same notation for its restriction to a Borel measurable subset of $\R^d$ with non-empty interior and boundary that is a null-set for $\lambda$. Accordingly, in a topological setting, measurability of sets and functions will refer to the respective Borel $\sigma$-algebras.

The space of measures $\CM(S)$ is equipped with the Fortet-Mourier norm $\|\mu\|_\FM^*$:
\begin{equation}\label{eq:def Fortet-Mourier norm}
    \|\mu\|_\FM^* := \sup\bigl\{ \int_S f\,d\mu\colon f\in\BL(S),\ \|f\|_\infty\leq 1,\ |f|_L\leq 1\bigr\}.
\end{equation}
It can be readily checked that $\|f\cdot\mu\|_\FM^*\leq 2\|f\|_\FM \|\mu\|_\FM^*$, for any $\mu\in\CM(S)$ and $f\in\BL(S)$.
On $\CM^+(S)$ the Fortet-Mourier norm metrizes  the usual weak topology on measures furnished by pairing $\CM(S)$ with $C_b(S)$ via integration:
\begin{equation}\label{eq:pairing}
    \langle\mu, f\rangle := \int_S f\,d\mu,\qquad \mu\in\CM(S),\ f\in C_b(S).
\end{equation}
$\CM^+(S)$ is closed in $\CM(S)$ for this weak topology. 
The metric derived from the Fortet-Mourier norm turns $\CM^+(S)$ into a complete separable metric space \cite[Section 1.3]{dull2021spaces}. 
We denote by $\CM(S)_\FM$ and $\CM^+(S)_\FM$ the respective spaces of measures equipped with the Fortet-Mourier metric and corresponding metric topology. It is well known that $\CM(S)_\FM$ is not complete \cite{gwiazda2018measures}, generally. 
Its metric completion is denoted by $\overline{\CM}(S)_\FM$, viewed as closure of $\CM(S)$ in $\BL(S)^*$, through the pairing \eqref{eq:pairing}. 

The total variation norm,
\[
    \|\mu\|_\TV\ :=\ |\mu|(S)\ =\ \mu^+(S) + \mu^-(S),\qquad \mu\in\CM(S),
\]
where $\mu^+$ and $\mu^-$ are the positive and negative part of the signed measure $\mu$ in the Hahn-Jordan decomposition, is closely related to the order structure on $\CM(S)$ -- defined set-wise -- and turns $\CM(S)$ into a Banach lattice (cf. e.g. \cite{Meyer-Nieberg:1991}, Section 1.1, Example {\it vi}, p.9). One has
\begin{equation}\label{eq:expression TV-norm}
    \|\mu\|_\TV\ =\ \sup\bigl\{ \langle\mu,f\rangle\colon f\in B_b(S),\ \|f\|_\infty\leq 1\bigr\}\ =\ \sup\bigl\{ \langle\mu,f\rangle\colon f\in C_b(S),\ \|f\|_\infty\leq 1\bigr\}. 
\end{equation}
Thus, \eqref{eq:def Fortet-Mourier norm} and \eqref{eq:expression TV-norm} combined yield that $\|\cdot\|_\TV$ is stronger than $\|\cdot\|_\FM^*$: one has $\|\mu\|_\FM^*\leq\|\mu\|_\TV$. It follows, that for $\mu\in\CM^+(S)$, $\|\mu\|_\TV=\|\mu\|_\FM^*$. For any $R>0$, we define
\[
    \CM_R(S) := \bigl\{ \mu\in\CM(S)\colon \|\mu\|_\TV\leq R\bigr\}.
\]
The relevance of this ball of radius $R$ in $\|\cdot\|_\TV$-norm in what is to follow lies in `repairing' to some extent the lack of completeness of $\CM(S)_\FM$. That is,
\begin{prop} \label{prop:completeness TV-ball}
    For any $R>0$, $\CM_R(S)$ is complete in $\bigl( \CM(S),\|\cdot\|_\FM^*\bigr)$.
\end{prop}

\begin{proof} 
    For $S=\R_+$, this result can be found in \cite[Theorem 2.7]{GWIAZDA20102703}, proven by a duality argument, which does not work for general Polish spaces. We shall prove the general case.
    Let $(\mu_n)_n\subset\CM_R(S)$ be a Cauchy sequence for $\|\cdot\|_\FM^*$. Let $\phi\in\overline{\CM}(S)$ be its limit. For every $f\in\BL(S)$, $\langle \mu_n,f\rangle\to\langle \phi,f\rangle$. $\BL(S)$ is $\|\cdot\|_\infty$-dense in the space $U_b(S)$ of bounded uniformly continuous functions on $S$ (cf. \cite{Dudley:1966}, Lemma 8). Because of $\sup_n\|\mu_n\|_\TV\leq R$ and \eqref{eq:expression TV-norm}, $(\langle \mu_n, g\rangle)_n$ is a Cauchy sequence in $\R$ for any $g\in U_b(S)$. Thus, $\bar{\phi}(g):= \lim_{n\to\infty} \langle \mu_n, g\rangle$ defines a linear functional on $U_b(S)$. The Banach-Steinhaus Theorem yields that $\bar{\phi}\in U_b(S)^*$. Clearly, the restriction of $\bar{\phi}$ to $\BL(S)$ equals $\phi$. Pachl has shown that $\CM(S)$ is $U_b(S)$-weakly sequentially complete (cf. \cite{Pachl:1979} Theorem 3.2, or \cite{Pachl:2013} Theorem 5.45). Hence, $\bar{\phi}$ is represented by integration against some $\mu\in\CM(S)$ -- and the same holds for $\phi$. It can be verified that
    \begin{equation}
        \|\mu\|_\TV\ =\ \sup\bigl\{ \langle\mu,g\rangle \colon g\in U_b(S),\ \|g\|_\infty\leq 1\bigr\}.
    \end{equation}
    Therefore, $\|\mu\|_\TV\leq R$ and we conclude that $\CM_R(S)$ is $\|\cdot\|_\FM^*$-complete.
\end{proof}

If $\Phi:S\to S'$ is a measurable map between measurable spaces, it defines a map between the spaces of signed measures and between the convex cones of positive measures on $S$ and $S'$ by push-forward. This map is denoted by $\Phi_\#$. Thus, for any positive or finite signed measure $\mu$ on $S$ , $\Phi_\#\mu$ is such a measure on $S'$, defined for any measurable set $A\subset S'$ by 
\[
    \Phi_\# \mu (A) := \mu(\Phi^{-1}(A)).
\]
We denote the (finite) positive measure on the torus $\mathbb{T}^d$, obtained by push-forward of $\lambda$ under the canonical quotient map, by $\lambda$ as well.

In the case that transformation $\Phi$ maps $S$ into itself, then $\Phi_\#$ defines a {\it Markov operator} on $\CM^+(S)$ or $\CM(S)$ which we shall denote by $P_\Phi$. That is, it is an operator that maps $\CM^+(S)$ into itself, is additive and positive-scalar homogeneous, and preserves total mass of positive measures.
The operator $P_\Phi$ is {\it regular}, with dual operator $U_\Phi f = f\circ\Phi$, meaning that
\[
    \langle P_\Phi\mu, f\rangle\ =\ \int_S f \,d(\Phi_\#\mu)\ =\ \int_S f\circ\Phi\,d\mu = \langle\mu, U_\Phi f\rangle\qquad \mbox{for all}\ \mu\in\CM(S),\ f\in B_b(S).
\]

We shall make use of Bochner integration of Banach-space-valued functions, see \cite{cohn2013measure, Diestel-Uhl} for details. A map $f$ from a positive measure space $(S,\Sigma,\mu)$ into a Banach space $X$ is strongly Bochner measurable if it is the point-wise limit limit of a sequence of simple functions $f_n$, i.e. of the form
\[
    f_n(s)\ =\ \sum_{k=1}^{N_n} x_k^{(n)}\, \ind_{E^{(n)}_k}(s),
\]
where $\ind_A$ is the indicator function of the $A$, ${E^{(n)}_k}\in\Sigma$ and $x_k^{(n)}\in X$. A function $f$ is Bochner integrable with respect to $\mu$ if and only if it is strongly Bochner measurable and $s\mapsto \|f(s)\|_X$ is $\mu$-integrable. Then,
\[
    \left\| \int_E f(s)\,d\mu(s)\right\|_X\ \leq\ \int_E \bigl\| f(s) \|_X\,d\mu(s).
\]
The following result is remarkable, because a similar `triangle inequality for integrals' holds for the total variation norm for Bochner integration in the Banach space $\overline{\CM}(S)_\FM$, for which it obviously is not the norm.
\begin{prop}[\cite{evers2015_JDE}, Proposition C.2] \label{prop:measure Bochner integral and TV-norm}
    Let $\nu$ be a $\sigma$-finite positive measure on the measurable space $(S',\Sigma')$ and let $f:S'\to \CM(S)_\FM$ be strongly Bochner measurable with respect to $\nu$ in $\overline{\CM}(S)_\FM$. Then $s'\mapsto \|f(s')\|_\TV$ is $\Sigma'$-Borel measurable. If it is $\nu$-integrable, then $f$ is Bochner integrable in $\overline{\CM}(S)_\FM$ with respect to $\nu$ and 
    \[
        \left\| \int_{S'} f(s')\,d\nu(s') \right\|_\TV\ \leq \int_{S'} \bigl\| f(s') \bigr\|_\TV\,d\nu(s').
    \]
    In particular, the Bochner integral of $f$ is a finite (signed) measure.
\end{prop}

\subsection{Fundamental assumptions and solution concept}
\label{subsec:solution operator}

The mathematical literature provides several well-posedness results for problem (P) or for variants of it, see for example \cite{ackleh2020well}, where a more complicated system of equations is studied with density dependent velocity field that is otherwise not dependent on time -- or \cite{Canizo_ea:2012} for the case of a time-independent and density independent velocity field. `Our' case of a time-dependent, density-independent velocity field is neither covered in \cite{ackleh2020well,Canizo_ea:2012},  nor in \cite{Ambrosio-Crippa:2014}, where there is no reaction term. 
 
 In particular, we will be following \cite{ackleh2020well} by assuming throughout this work that the velocity field satisfies:
\begin{quote}
\begin{enumerate}
    \item[(V1)] {\it $(t,x)\mapsto v_t(x)\in \R^d$ is continuous;}
    \item[(V2)] {\it $x\mapsto v_t(x)$ is bounded and globally Lipschitz on $\Omega$ for each $t\geq 0$;}
    \item[(V3)] {\it$t\mapsto v_t: \R_+\to \BL(\Omega,\R^d)$ is continuous.}
\end{enumerate}
\end{quote}
The Lipschitz constant of $v_t$ is $|v_t|_L$. As norm on $\BL(\Omega,\R^d)$ we take $\|v\|_\FM := \max(\|v\|_\infty,|v|_L)$. We require (V3) to ensure that for any $T>0$, there exists $L_T>0$ such that $|v_t|_L\leq L_T$ for all $t\in[0,T]$. This is a condition imposed in \cite{ackleh2020well} in their results on well-posedness for problem (P) (see e.g. \cite{Ambrosio-Crippa:2014}, where the more general assumption that $v\in L^1(\R_+,\BL(\Omega,\R^d))$ is made).

With conditions (V1)--(V3), $v$ generates a flow $\Phi^v_{s,t}$ on $\Omega$ for $s,t\in\R$ with $0\leq s\leq t$. That is, the motion $t\mapsto\Phi_{s,t}^v(x_0)$ is the unique and globally existing solution to 
\begin{equation}\label{Eq:FlowODE}
    \left\lbrace
    \begin{split}
        \Dot{x}(t) &= v_t(x(t))\quad \mbox{for}\ t>s,\\
        x(s) &=x_0\in\Omega
    \end{split}. \right.
\end{equation}
The flow maps $\Phi^v_{s,t}$ satisfy $\Phi^v_{t,t'}\circ\Phi^v_{s,t}=\Phi^v_{s,t'}$ for $t'\geq t\geq s\geq 0$. Moreover, they are Lipschitz, with
\begin{equation}\label{eq:estimate Lipschitz const flow}
    \bigl| \Phi_{s,t}^v \bigr|_L\ \leq\ \exp\bigl( \int_s^t |v_\tau|_L\,d\tau\bigr)\  =:\ L^v_{s,t}
\end{equation}
(cf. \cite{Ambrosio-Crippa:2014}, Eq. (2.1)).

We further associate to $v_\bullet$ the following family of Markov operators on measures by means of push-forward under the flow:
\begin{equation}
    P_{s,t}^v \mu\ :=  (\Phi_{s,t}^v)_\# \mu,\qquad 0\leq s\leq t.
\end{equation}
The linear operators $P^v_{t_0,t}$ are continuous on $\CM(S)_\FM$. They satisfy the norm estimate
\begin{equation}\label{eq:estimate FM-operator norm MOflow} 
    \bigl\| P^v_{t_0,t}\mu \bigr\|_\FM^*\ \leq\ \|\mu\|_\FM^*\,\max\bigl(1, |\Phi^v_{t_0,t}|_L\bigr)\ \leq\ L^v_{t_0,t}\|\mu\|_\FM^* \qquad \mbox{for all}\ \mu\in\CM(\Omega).
\end{equation}
It is clear from their definition that the operators $P_{s,t}^v$ preserve $\CM^+(\Omega)$.

It is well known that $\mu_{t}:= P^v_{t_0,t}\nu$ solves the linear continuity equation 
\begin{equation}\label{eq:linear conservation law}
    \partial_t\mu_t + \mathrm{div}_x\bigl( v_t(x)\mu_t\bigr) =0,\qquad \mu_{t_0}=\nu,
\end{equation}
(see e.g. \cite{villani2003topics} Theorem 5.34, p.167). 
The operators $P^v_{s,t}$ are therefore the building blocks for representing the solutions to Problem (P). First one considers the perturbation of \eqref{eq:linear conservation law}, by introducing a non-zero right-hand side, independent on the density $\mu_t$:
\begin{equation}\label{eq:linear transport with forcing}
    \partial_t\mu_t + \mathrm{div}_x\bigl( v_t(x)\mu_t\bigr) = \sigma_t,\qquad \mu_{t_0}=\nu,
\end{equation}
where $\sigma_t\in\CM(\Omega)$. One has the following well-posedness result:
\begin{prop}[{\bf Ackleh \& Saintier \cite{ackleh2020well}, Proposition 5.1}]\label{prop:existence solution forced eq}
    Let $v_\bullet$ satisfy (V1)--(V3). let $T>0$ and assume that $\sigma_\bullet\in C([0,T], \CM(\Omega)_\FM)$ is such that $\sup_{0\leq t\leq T} \|\sigma_t\|_\TV<\infty$. Then, for any $t_0\in[0,T)$ and initial condition $\nu\in\CM(\Omega)$, equation \eqref{eq:linear transport with forcing} has a unique solution $\mu_\bullet\in C([t_0,T],\CM(\Omega)_\FM)$. It is given explicitly by
    \begin{equation}\label{eq:VoC formula}
        \mu_t\ = P^v_{t_0,t}\nu\ +\ \int_{t_0}^t P^v_{s,t}\sigma_s\,ds.
    \end{equation}
\end{prop}
\noindent The integral in \eqref{eq:VoC formula} should be interpreted as a Bochner integral in $\overline{\CM}(\Omega)_\FM$, which -- as part of the result -- takes values in $\CM(\Omega)$.
\vskip 2mm

We wish to consider problem (P) for any initial time $t_0\in [0,T)$, i.e. require $\mu_{t_0}=\nu$. Accordingly,
\begin{definition}
    A {\bf mild solution} to problem (P) on $[t_0,T]$ for the initial condition $\mu_{t_0}=\nu$ is a continuous map $\mu_\bullet:[t_0,T]\to \CM(\Omega)_\FM$ that satisfies 
    \begin{equation}\label{eq:VoC formula mild solution}
        \mu_t = P^v_{t_0,t}\nu + \int_{t_0}^t P^v_{s,t}[f_s(\mu_s)]\,ds,\qquad t\in[t_0,T],
    \end{equation}
    and that is bounded in total variation norm, i.e. $\sup_{t\in[t_0,T]}\|\mu_t\|_\TV<\infty$.
\end{definition}
The integral in \eqref{eq:VoC formula mild solution} must be interpreted as Bochner integral as above, as for Proposition \ref{prop:existence solution forced eq}. For a continuous map $t\mapsto \mu_t$ into $\CM(\Omega)_\FM$, the function $t\mapsto \|\mu_t\|_\TV$ is Borel measurable (cf. Lemma \ref{lem:Jordan is measurable} and Lemma \ref{lem:measurability postive part}). There is no {\it a priori} guarantee, that it is bounded. Therefore, it is included in the definition. If $\mu_t$ is positive for every $t$, then continuity of $\mu_\bullet$ assures bounded total variation, since $\|\mu_t\|_\FM^*=\|\mu_t\|_\TV$ for $\mu_t\in\CM^+(\Omega)$.

The integrand in \eqref{eq:VoC formula mild solution} -- as a measure-valued function of $s$ -- must be strongly Bochner measurable. The conditions imposed on $f_\bullet$, following \cite{Canizo_ea:2012,ackleh2020well}, will ensure this:
\begin{quote}
    \begin{enumerate}
        \item[(A1)] {\it $f:(t,\mu)\mapsto f_t(\mu):\R_+\times \CM(\Omega)_\FM\to\CM(\Omega)_\FM$ is jointly continuous;}
        \item[(A2)] {\it For every $t\geq0$ and $R>0$ there exists $L_R^f>0$ such that 
        \[
            \bigl\| f_t(\mu_1)-f_t(\mu_2) \bigr\|_\FM^* \leq L_R^f\|\mu_1-\mu_2\|_\FM^*\qquad \mbox{for all}\ \mu_i\in\CM_R(\Omega).
        \]}
        \item[(A3)] {\it For every $R>0$ there exists $C_R^f>0$ such that 
        \[
           \|f_t(\mu)\|_\TV \leq C_R^f \qquad \mbox{for all}\ \mu\in\CM_R(\Omega),\ t\geq 0.
        \]}
    \end{enumerate}
\end{quote}
\noindent In fact, (A1) and (A2) imply that $s\mapsto f_s(\mu_s)$ is continuous for continuous $\mu_\bullet$ and  that consequently the Bochner integral in \eqref{eq:VoC formula mild solution} is well-defined and exists. Note that the function $R\mapsto C^f_R$ in (A3) must be non-decreasing. It is immediately clear that any finite linear combination of maps that satisfy (A1)--(A3) will satisfy the very same assumptions.
\vskip 2mm

Examples of reaction maps $f_\bullet$ that satisfy Assumptions (A1)--(A3) appear naturally in many modelling settings, as the following proposition illustrates. The map $f_\bullet$ is split there in two parts: a positive general production term (`$p$') and a term that combines degradation and individual-induced proliferation. The latter should then be absolutely continuous with respect to the given mass distribution, since by the latter processes mass can only appear or disappear where there is mass. 
\begin{prop}\label{prop:example family admissible reaction terms}
    Let $p_\bullet:\R_+\times \CM(\Omega)_\FM\to\CM^+(\Omega)_\FM:(t,\mu)\mapsto p_t(\mu)$ satisfy (A1)--(A3). Let $F_\bullet:\R_+\times\CM(\Omega)_\FM\to \BL(\Omega):(t,\mu)\mapsto F_t(\mu)$ be jointly continuous and satisfies:
    \begin{enumerate}
        \item[(A2')] For every $t\geq 0$ and $R>0$ there exists $\hat{L}^F_R>0$ such that
        \[
            \bigl\|F_t(\mu_1)-F_t(\mu_2) \bigr\|_\FM\ \leq\ \hat{L}_R^F \|\mu_1-\mu_2\|_\FM^*\qquad \mbox{for all}\ \mu_i\in\CM_R(\Omega);
        \]
        \item[(A3')] For every $R>0$ there exists $\hat{C}^F_R>0$ such that
        \[
            \|F_t(\mu)\|_\infty\ \leq\ \hat{C}_R^F\qquad \mbox{for all}\ \mu\in\CM_R(\Omega),\ t\geq 0.
        \]
    \end{enumerate}  
    \noindent Then, 
    \begin{equation}\label{eq:specific admiss non-linearity}
        f_t(\mu)\ :=\ p_t(\mu)\ +\ F_t(\mu)\cdot \mu
    \end{equation}
    satisfies (A1)--(A3).
\end{prop}
\begin{proof}
     Recall that $\|g\cdot\mu\|_\FM^*\leq 2\|g\|_\FM\|\mu\|_\FM^*$ for any $g\in\BL(\Omega)$ and $\mu\in\CM(\Omega)$. Hence, for $g_1,g_2\in\BL(\Omega)$ and $\mu_1,\mu_2\in\CM(\Omega)$,
     \[
        \bigl\| g_1\cdot\mu_1 - g_2\cdot\mu_2 \bigr\|_\FM^* \ \leq\ 2\|g_1-g_2\|_\FM \|\mu_1\|_\FM^*\ + \ 2\|g_2\|_\FM\|\mu_1-\mu_2\|_\FM^*.
     \]
     Thus, the bilinear map $(g,\mu)\mapsto g\cdot \mu:\BL(\Omega)\times\CM(\Omega)_\FM\to\CM(\Omega)_\FM$ is jointly continuous.  Then, for $\mu_1,\mu_2\in\CM_R(\Omega)$ and $t\geq 0$,
     \[
        \bigl\| F_t(\mu_1)\cdot\mu_1 - F_t(\mu_2)\cdot\mu_2 \bigr\|_\FM^*\  \leq\ (2R \hat{L}_R^F + 2\hat{C}_R^F)\, \|\mu_1-\mu_2\|_\FM^*,
     \]
     so (A2) is satisfied. Moreover, for $\mu\in\CM_R(\Omega)$ and $t\geq 0$,
     \[
        \bigl\| F_t(\mu)\cdot\mu \bigr\|_\TV\ \leq \bigl\| F_t(\mu) \bigr\|_\infty\, \|\mu\|_\TV\ \leq\ R\hat{C}_R^F,
     \]
     which proves (A3).
\end{proof}

\subsection{Well-posedness -- Revisiting the construction of the solution operator}

The well-posedness of Problem (P) or, equivalently, of Problem \eqref{eq:VoC formula mild solution} can be shown along the lines of reasoning of Section 6 in \cite{ackleh2020well} or the proof of Theorem 2.4 in \cite{Canizo_ea:2012}. We will get into some details here, because some intermediate results or constructions will be of relevance later. 

\begin{lemma}[{\it Uniqueness}] \label{lem:uniqueness mild solutions}
    Assume that $\mu^1_\bullet$ and $\mu^2_\bullet$ are both mild solutions to Problem (P) on $[t_0,T]$ for initial condition $\nu\in\CM(\Omega)$. Then $\mu^1_t=\mu^2_t$ for all $t\in[t_0,T]$.
\end{lemma}
\begin{proof}
    Let $R:=\max_i\, \sup_{t\in[t_0,T]} \|\mu^i_t\|_\TV<\infty$. Then \eqref{eq:VoC formula mild solution} yields for any $t\in [t_0,T]$:
    \begin{align*}
        \|\mu^1_t - \mu^2_t\|^*_\FM \ &\leq\ \int_{t_0}^t \left \| P^v_{s,t}\bigl[ f_s(\mu^1_s) - f_s(\mu^2_s)\bigr] \right\|_\FM^*\, ds\ 
        \leq \ \int_{t_0}^t  L^v_{s,t} \bigl\| f_s(\mu^1_s) - f_s(\mu^2_s) \bigr\|_\FM^*\,ds\\
        &\leq \ L_R^f\int_{t_0}^t L^v_{s,t} \bigl\| \mu^1_s - \mu^2_s \bigr\|_\FM^*\,ds
    \end{align*}
    Application of a Gr\"onwall's Lemma gives $\|\mu^1_t - \mu^2_t\|^*_\FM=0$ for $t\in [t_0,T]$.
\end{proof}

For fixed vector field $v_\bullet$, $\delta>0$ and initial condition $\nu\in\CM(\Omega)$, let $t_0\geq 0$ and $\tau>0$ and put $I_\tau:=[t_0,t_0+\tau]$. The argumentation consists essentially of considering the iteration of operators $\Psi_{t_0,\tau}: C_{\nu,\delta}(I_\tau,\CM(\Omega)_\FM) \to C_{\nu,\delta}(I_\tau,\CM(\Omega)_\FM)$, where
\[
    C_{\nu,\delta}(I_\tau,\CM(\Omega)_\FM)\ :=\ \bigl\{ \mu_\bullet\in C(I_\tau,\CM(\Omega)_\FM)\colon \mu_{t_0}=\nu,\ \sup_{t\in I_\tau}\|\mu_t\|_\TV\leq \|\nu\|_\TV+\delta \bigr\}
\]
comes with the the metric associated to the supremum norm on $C(I_\tau,\CM(\Omega)_\FM)$. In view of Proposition \ref{prop:completeness TV-ball}, the metric space $C_{\nu,\delta}(I_\tau,\CM(\Omega)_\FM)$ is complete. $\Psi_{t_0,\tau}$ is defined by
\begin{equation}\label{eq:Picard iteration operator}
    \Psi_{t_0,\tau}(\mu_\bullet)(t)\ :=\ P^v_{t_0,t}\mu_{t_0} + \int_{t_0}^t P^v_{s,t}\bigl[f_s(\mu_s)\bigr]\,ds,\qquad t\in I_\tau.
\end{equation}
The conditions (V1)--(V3) on $v_\bullet$ and (A1)--(A3) on $f_\bullet$ ensure that for sufficiently small $\tau$, depending only on $\nu$ through $\|\nu\|_\TV$, $\Psi_{t_0,\tau}$ is a strict contraction on $C_{\nu,\delta}(I_\tau,\CM(\Omega)_\FM)$. The Banach Fixed Point Theorem then yields the existence of a (unique) mild solution on $I_\tau$. We shall sketch details of the argument through a few main observations and estimates.
\vskip 2mm

The first auxiliary result is:
\begin{lemma}
    Assume that $v_\bullet$ and $f_\bullet$ satisfy (V1)--(V3) and (A1)--(A3). Let $t_0\geq 0$, $R>0$ and $\nu\in\CM(\Omega)$ and put $r:=\|\nu\|_\TV+\delta$. Then the following statements hold:
    \begin{enumerate}
        \item[(i)] If $\mu_\bullet\in C_{\nu,\delta}(I_\tau,\CM(\Omega)_\FM)$, then $\|\Psi_{t_0,\tau}(\mu_\bullet)(t)\|_\TV \leq \|\nu\|_\TV + C_r^f(t-t_0)$ for any $t\in I_\tau.$
        \item[(ii)] If $\tau\leq \delta / C^f_{R+\delta}$, then $\Psi_{t_0,\tau}$ maps $C_{\nu,\delta}(I_\tau,\CM(\Omega)_\FM)$ into itself, for any $\nu\in\CM_R(\Omega)$.
    \end{enumerate}
\end{lemma}
\begin{proof}
    {\it (i).}\ One has for $t\in[t_0,T]$, using Proposition \ref{prop:measure Bochner integral and TV-norm}:
    \begin{align*}
        \|\Psi_{t_0,\tau}(\mu_\bullet)(t)\|_\TV \ &\leq\ \| P^v_{t_0,t}\nu\|_\TV\ +\ \left\| \int_{t_0}^t P^v_{s,t}\bigl[ f_s(\mu_s)\bigr]\,ds \right\|_\TV
        \ \leq\ \|\nu\|_\TV\ +\ \int_{t_0}^t \left\| P^v_{s,t}\bigl[ f_s(\mu_s)\bigr]\right\|_\TV\,ds\\
        &\leq \ \|\nu\|_\TV\ +\ \int_{t_0}^t \left\| f_s(\mu_s)\right\|_\TV\,ds\ 
        \leq \ \|\nu\|_\TV\ +\ C_r^f(t-t_0).
    \end{align*}
    {\it (ii)} is immediately clear from {\it (i)} and $x\mapsto C^f_x$ being non-decreasing on $(0,\infty)$.
\end{proof}

The existence of a (necessarily unique) fixed point for $\Psi_{t_0,\tau}$ in $C_{\nu,\delta}(I_\tau,\CM(\Omega)_\FM)$ -- hence, a mild solution to Problem (P) on $I_\tau$ -- for any $\nu\in\CM_R(\Omega)$ is guaranteed by making $\tau$ smaller than $\delta/C^f_{R+\delta}$, in such a way that $\Psi_{t_0,\tau}$ becomes a strict contraction on each of the closed, hence complete, sets $C_{\nu,\delta}(I_\tau,\CM(\Omega)_\FM)$, $\nu\in\CM_R(\Omega)$ (cf. \cite{ackleh2020well} Section 6). Thus, for every $\nu\in\CM_R(\Omega)$ one obtains a sequence $(\tau_n)_{n\in\NN}$ with $0<\tau_n\leq \delta C^f_{R+n\delta}$, $t_n:= t_{n-1}+\tau_n$, and a `chain' of mild solutions $\mu^n_\bullet$ to Problem (P) on $[t_{n-1},t_n]$  with initial condition $\mu^n_{t_{n-1}} = \mu^{n-1}_{t_{n-1}}$, where $\mu^1_{t_0}=\nu\in \CM_R(\Omega)$. Each $\mu^n_\bullet$ is the fixed point of $\Psi_{t_{n-1},\tau_n}$. One has
\[
    \sup_{t_{n-1}\leq t \leq t_n} \bigl\| \mu^n_t \bigr\|_\TV\ \leq R + n\delta.
\]
Piecing the parts of the chain continuously together, one obtains a mild solution $\mu_\bullet$ to Problem (P) with initial condition $\nu$ at $t_0$, on the interval $[t_0, \tilde{T})$, where $\tilde{T}:= t_0 + \sum_{n=1}^\infty \tau_n$. If $\tilde{\nu}:=\lim_{t\to \tilde{T}} \mu_t$ exists in $\overline{\CM}(\Omega)_\FM$, then $\|\tilde{\nu}\|_\FM^*<\infty$. If $\|\tilde{\nu}\|_\TV<\infty$ too, then the mild solution can be extended from initial (measure) condition $\tilde{\nu}$ at time $\tilde{T}$, following a similar procedure. Otherwise, $\tilde{\nu}$ is not a signed measure, and the mild solution cannot be extended beyond $\tilde{T}$.

Thus, putting all pieces together, one arrives at the following result:

\begin{prop}[{\bf Well-posedness}; \cite{ackleh2020well} Section 6 and \cite{Canizo_ea:2012}, Theorem 2.4]\,\\
    Assume that $v_\bullet$ and $f_\bullet$ satisfy (V1)--(V3) and (A1)--(A3). Let $t_0\geq 0$ and $\nu\in\CM(\Omega)$ a given initial condition at time $t_0$. There exists a maximal time $T^*>t_0$, depending on $\|\nu\|_\TV$ and $t_0$, such that there is a unique mild solution $\mu_\bullet\in C([t_0,T^*),\CM(\Omega)_\FM)$ to Problem (P) that is of bounded total variation on any $[t_0,T]$ for $t_0<T<T^*$. One has:
    \begin{enumerate}
        \item[(i)] Either $T^*=\infty$, or $\lim_{t\uparrow T^*} \|\mu_t\|_\TV = +\infty$.
        \item[(ii)] $\mu_\bullet$ depends locally Lipschitz-continuously on $\nu$ in $\|\cdot\|^*_\FM$-norm, uniformly locally in time. That is, for every $T\in (t_0,T^*)$ there exist $R>0$ and there exists $\omega_{R,T}\geq 0$ such that for mild solutions $\mu^i_\bullet$ with initial conditions $\nu_i$ that satisfy $\|\nu_i\|_\TV\leq R$, $i=1,2$, one has
        \[
            \|\mu^1_t - \mu^2_t\|_\FM^*\ \leq\ e^{\omega_{R,T} t}
            \,\|\nu_1-\nu_2\|_\FM^*,\qquad\mbox{for all}\ t\in[t_0,T].
        \]
    \end{enumerate}
\end{prop}
\begin{proof}
    {\it (i).}\ Standard argument by contradiction. If $T^*<\infty$, and the limit is not infinite, one can extend the solution slightly beyond $T^*$; see e.g. the extension arguments used in \cite{Amann1990} when identifying global solutions to differential equations.\\
    {\it (ii).} For any $t\in[t_0,T^*)$, \eqref{eq:VoC formula mild solution} holds. Let $t_0<T<T^*$. Then, according to part {\it (i)}, there exists $R>0$, finite, such that $\|\mu^i_t\|_\TV\leq R$ for all $t\in[t_0,T]$, $i=1,2$. Then, one obtains, by applying  \eqref{eq:estimate Lipschitz const flow} and \eqref{eq:estimate FM-operator norm MOflow} in the following estimates:
    \begin{align*}
        \|\mu^1_t - \mu^2_t\|^*_\FM \ &\leq\ \bigl\| P^v_{t_0,t}(\nu_1-\nu_2)\bigr\|_\FM^*\ +\ \int_{t_0}^t \left \| P^v_{s,t}\bigl[ f_s(\mu^1_s) - f_s(\mu^2_s)\bigr] \right\|_\FM^*\, ds\\
        &\leq \ L^v_{t_0,t} \|\nu_1-\nu_2\|_\FM^* + \int_{t_0}^t L^v_{s,t} \bigl\| f_s(\mu^1_s) - f_s(\mu^2_s) \bigr\|_\FM^*\,ds\\
        &\leq \ L^v_{t_0,t}\|\nu_1-\nu_2\|_\FM^*+ L_R^f\int_{t_0}^t L^v_{s,t} \bigl\| \mu^1_s - \mu^2_s \bigr\|_\FM^*\,ds
    \end{align*}
    Application of a Grönwall's Inequality yields the result. 
\end{proof}

\begin{remark}\label{rem:solution locally Picard iteration}
Thus, let $v_\bullet$ satisfy (V1)--(V3) and let $f_\bullet$ satisfy (A1)--(A3). Take $t_0\geq 0$ and $R>0$. The maximal mild solution for initial condition $\nu\in\CM_R(\Omega)$ is obtained by piecing together local mild solutions on intervals $I_{t',\tau}:= [t',t'+\tau]$ in the maximal domain of existence. On such an interval, the solution is the limit of a Picard iteration procedure. Thus, for any $t'$ in the maximal interval of existence, there exists $\tau$ such that the restriction of $\mu_\bullet$ to $I_{t',\tau}$ is the unique fixed point of $\Psi_{t',\tau}$ in $C_{\nu',\delta}(I_{t',\tau},\CM(\Omega)_\FM)$, obtained by iteration of the operator $\Psi_{t',\tau}$:
\begin{equation}\label{eq:Picard iteration in measures}
    \mu_\bullet = \lim_{n\to \infty} \mu^{(n)}_\bullet,\qquad\quad \mu^{(0)}_t := P^v_{t',t}\nu'\ \mbox{for}\ t'\leq t\leq t'+\tau,\quad \mu^{(n+1)}_\bullet:= \Psi_{t',\tau}\bigl(\mu^{(n)}_\bullet\bigr).
\end{equation}
Note that the limit is that of uniform convergence on $I_{t',\tau}$ in $\|\cdot\|_\FM^*$-norm in $\CM(\Omega)$.
\end{remark}

\section{Weak limit of sequences of measures with $L^p$-densities}\label{weaks}


In this section, we collect results that hold in the general setting of Polish spaces, which we denote by $S$ to distinguish them from the specific domain $\Omega$. These are novel and of general interest, we believe, because there are hardly any results in the literature that connect weak convergence in measures -- for measures with $L^p$-densities -- to absolute continuity of the limit measure, with similar regularity for its density.

The main result will precisely allow us to conclude this -- under suitable conditions of course: the limit in the $C_b$-weak topology on measures of a sequence of measures, each absolutely continuous with respect to the same reference measure and with Radon-Nikodym derivative in $L^p$, is itself in $L^p$ too, provided $1<p\leq\infty$.  If $p=1$, then  such a result cannot be expected to hold. Take for example on $\R$ the sequence of probability measures that are normally distributed with mean 0 and with standard deviation converging to zero. Then, this sequence of measures converges weakly to $\delta_0$, while they all have $L^1$-norm $1$.

The proof of our result is founded on a characterisation of $L^p$-functions encountered in \cite{BogachevI:2007} Ex. 4.7.102, p.315:

\begin{lemma}
    Let $\nu$ be a finite positive measure on a measurable space $(S,\Sigma)$. Then, a function $f$ is in $L^p(S,\nu)$ with $1<p<\infty$ if and only if there exists $C>0$ such that for every finite disjoint measurable partition $A_1,\dots, A_m$ of $S$ with $\nu(A_k)>0$ for all $k$,
    \begin{equation}
        \sum_{k=1}^m \nu(A_k)^{1-p} \left| \int_{A_k} f d\nu \right|^p\ \leq \ C.
    \end{equation}
    The smallest such constant $C$ is $\|f\|_p^p$.
\end{lemma}

Consider now a general Polish space $S$, equipped with its Borel $\sigma$-algebra. Let $\nu$ be a positive Borel measure on $S$. Define
\[
    L^p_+(\nu) := \bigl\{ f\in L^1\cap L^p(S,\nu)\colon f\geq 0\ \mbox{$\nu$-a.e} \bigr\}.
\]
The crucial auxiliary result alluded to above is the following.
\vskip 2mm

\begin{theorem}\label{Lem:BoundedWeakLimits}
    Let $\nu$ be a positive Radon measure on the Polish space $S$. Let $1<p\leq \infty$ and $f_n\in L^p_+(\nu)$ for $n\in\NN$. Let $\mu_n:= f_nd\nu\in\CM^+(S)$. Suppose that $\mu_n\to\mu_*$ in the weak topology on measures.
    Then the following holds:
    \begin{enumerate}
        \item[(i)] There exists a $\sigma$-compact set $S_*$ on which all $\mu_n$ and $\mu_*$ are concentrated. The Borel measure $\nu_*$ defined by $\nu_*(E):=\nu(E\cap S_*)$ for any Borel $E\subset S$ is $\sigma$-finite.
        \item[(ii)] If $M:=\liminf_{n\to\infty}\|f_n\|_p<\infty$, then $\mu_*$ is absolutely continuous with respect to $\nu_*$ and hence $\nu$. The Radon-Nikodym derivative $f_*$ of $\mu_*$ with respect to $\nu_*$ is in $L^p(S,\nu_*)$, with $\|f_*\|_p\leq M$. In particular, we have $\mu_*=fd\nu$ with $f:=f_*\ind_{S_*}\in L^p_+(\nu)$ and $\|f\|_{L^p(\nu)}\leq M$.
    \end{enumerate}
\end{theorem}

\begin{proof} 
    {\it (i)}.\ Since $\mu_n\to\mu_*$ in the weak topology of measures, the set of measures $\CM_*:=\{\mu_n:n\in\NN\}\cup\{\mu_*\}$ is weakly (sequentially) compact. Hence, it is uniformly tight, by the Prokhorov Theorem (cf. \cite{BogachevII:2007}, Theorem 8.6.2, p.202). Therefore, there exists an increasing sequence $(K_m)_{m\in\NN}$ of compact sets in $S$, such that $\mu(S\setminus K_m)<\frac{1}{m}$ for all $\mu\in\CM_*$ and $m\in\NN$. Thus, all $\mu\in\CM_*$ are concentrated on the $\sigma$-compact set $S_*:=\bigcup_m K_m$ in $S$. Since $\nu$ is a Radon measure, $\nu(K_m)<\infty$. Without loss of generality, we may assume that $\nu(K_m)>0$ for all $m$.  In the usual way, one constructs from $(K_m)_{m\in\NN}$ a countable family $(E_k)_{k\in\NN}$ of pairwise disjoint Borel sets such that $K_m = \bigcup_{k=1}^m E_k$. Let $\nu_*$ be the Borel measure on $S$ obtained by restricting $\nu$ to $S_*$ as indicated. Define $E_0:=S\setminus S_*$. Then $\nu_*(E_m)<\infty$ for all $m\in\NN_0$, while $S=\bigcup_{m=0}^\infty E_m$. In particular, $\nu_*$ is $\sigma$-finite.
    
    \noindent{\it (ii)}.\ Since $\mu_n$ and $\mu_*$ are positive finite measures, Alexandrov's Theorem (cf. e.g. \cite{BogachevII:2007} Theorem 8.2.3, p.184) yields that for any $U\subset S$ open, 
    
    \begin{align}
        \mu_*(U)\ &\leq\ \liminf_{n\to\infty} \mu_n(U) = \liminf_{n\to\infty} \int_{U\cap S_*} f_n\,d\nu = \liminf_{n\to\infty} \sum_{m=1}^\infty \int_{U\cap E_m} f_n\,d\nu\nonumber\\
        &\leq\ \liminf_{n\to\infty} \sum_{m=1}^\infty \|f_n\|_p\, \nu(U\cap E_m)^{1/q}\ \leq\ M\nu_*(U)^{1/q},\label{eq:estimat mu by nu}
    \end{align}
    where we used Jensen's Inequality in the last step. In this context, $q$ is the coefficient conjugate to $p$: $\frac{1}{p}+\frac{1}{q}=1$. Since $S$ is Polish, $\mu_*$ is Radon and so is $\nu$ by assumption. Thus, both $\mu_*$ and $\nu_*$ are open outer regular on Borel sets $E$:
    \[
        \mu_*(E)\ =\ \inf\bigl\{ \mu_*(U)\colon U\subset S\ \mbox{open},\ E\subset U\bigl\}.
    \]
    Consequently, estimate \eqref{eq:estimat mu by nu} yields
    \begin{equation}\label{eq:estimate mu and nu}
        \mu_*(E) \ \leq\ M\nu_*(E)^{1/q},\quad \mbox{for all Borel}\ E\subset S.
    \end{equation}
    Hence, if $E\subset S$ is a $\nu_*$-null set, then $\mu_*(E)=0$. That is, $\mu_*$ is absolutely continuous with respect to $\nu_*$, $\mu_*\ll\nu_*$. If $E\subset S$ is a $\nu$-null set, then it is also a $\nu_*$-null set. So, also $\mu_*\ll\nu$.
    
    Since $\nu_*$ is $\sigma$-finite, the Radon-Nikodym Theorem implies the existence of the Radon-Nikodym derivative, $f_*\in L^1(S,\nu_*)$, of $\mu_*$ relative to $\nu_*$. 
    For the proof of the statements that $f_*\in L^p(S,\nu_*)$ and $\|f_*\|_p\leq M$, we need to distinguish the case $p=\infty$ from $1<p<\infty$.
    \vskip 2mm

    \noindent {\it Case $p=\infty$.}\ Letting $f:=f_*\ind_{S_*}$, one has
    \[
        \|f\|_\infty = \|f\|_{L^\infty(\nu_*)}\ =\ \inf_{S':\nu_*(S\setminus S')=0}\ \sup_{x\in S'} |f(x)|\ =\ \inf\bigl\{ M'\geq 0\colon |f|\leq M'\ \nu_*-\mbox{a.e.}\bigr\}.
    \]
    Suppose that $f_*$ is not essentially bounded. Then there must exists a Borel set $E$ and $m_0\in\NN$ such that $f_*\ind_E> M\ind_E$ and $\nu_*(E\cap E_{m_0})>0$. Then
    \[
        \mu_*(E\cap E_{m_0}) = \int_S f_*\ind_{E\cap E_{m_0}}\, d\nu_*\ >\ \int_S M\ind_{E\cap E_{m_0}}\, d\nu_* = M\nu_*(E\cap E_{m_0}), 
    \]
    which contradicts estimate \eqref{eq:estimate mu and nu} (note that $q=1$). So, $f_*$ must be essentially bounded. Similar argumentation as above, shows that the assumption that $\|f_*\|_\infty>M$ leads to a contradiction with \eqref{eq:estimate mu and nu}. So, the bound $\|f\|_\infty\leq M$ holds. 
    \vskip 2mm
    
    \noindent {\it Case $1<p<\infty$.}\ Let $q$ be the coefficient conjugate to $p$ and let $\nu_m$ be the Borel measure on $K_m$ obtained by restricting $\nu$ to the Borel $\sigma$-algebra $\mathcal{B}(K_m)$. We shall show that $f_*|_{K_m}\in L^p(K_m,\nu_m)$ and $\|f_*|_{K_m}\|_p\leq M$. Then, by monotone convergence,
    \[
        \int_S |f_*|^p\,d\nu_*\ =\ \lim_{m\to\infty} \int_{K_m} \bigl|f_*|_{K_m}\bigr|^p\,d\nu_m\ \leq\ M^p.
    \]
    To this end, fix $m\in\NN$ and let $A_1,\dots, A_N$ be a disjoint Borel measurable partitioning of $K_m$ with $\nu(A_k)>0$ for all $k$. 
    The assumptions on $f_n$ imply that
    \[
        \sum_{k=1}^N \nu(A_k)^{1-p} \left| \int_{A_k} f_nd\nu\right|^p \ \leq\ \|f_n\ind_{K_m}\|_{L^p(\nu_m)}^p\ \leq\ \|f_n\|_{L^p(\nu)}^p\quad\mbox{for all}\ n\in\NN.
    \]
    Let $\eps>0$. Since $\nu(A_k)<\infty$, for any Borel set $B$ in $S$ with $B\supset A_k$ and $n\in\NN$,
    \[
        0\leq \int_{B} f_nd\nu - \int_{A_k} f_nd\nu = \int_{B\setminus A_k} f_n d\nu\ \leq \ \|f_n\|_p\,\nu(B\setminus A_k)^{1/q}.
    \]
    Moreover, because $\nu$ is open outer regular, for each fixed $k$ one can select an open set $U_k$ in $S$ with $U_k\supset A_k$ and $\nu(U_k)<\infty$, such that
    \[
        \nu(U_k\setminus A_k)^{p/q}\ <\ \frac{1}{N}\nu(A_k)^{p/q}\cdot \frac{\eps^p}{M^p+1}.
    \]
    Then, applying Alexandrov's Theorem again to the measures $\mu_n = f_nd\nu\in\CM^+(S)$ that converge weakly to $\mu_* = f_* d\nu_*$,
    \begin{align*}
        \left( \sum_{k=1}^N \nu(A_k)^{1-p}\right. & \left.\left|\int_{U_k} f_* d\nu_*\right|^p \right)^{1/p}\ \leq\ \left( \sum_{k=1}^N \nu(A_k)^{1-p}\, \liminf_{n\to\infty} \left|\int_{U_k} f_n d\nu\right|^p \right)^{1/p}\\
        & \leq\ \liminf_{n\to\infty} \left( \sum_{k=1}^N \nu(A_k)^{1-p}\, \left|\int_{U_k} f_n d\nu\right|^p \right)^{1/p}\\
        & =\ \liminf_{n\to\infty}\left( \sum_{k=1}^N \left|\; \nu(A_k)^{-1/q} \int_{A_k} f_n d\nu\ +\ \nu(A_k)^{-1/q}\int_{U_k\setminus A_k} f_n d\nu \; \right|^p \right)^{1/p}\\
        & =\ \liminf_{n\to\infty} \left[\, \left(\sum_{k=1}^N \nu(A_k)^{1-p} \left|\int_{A_k}f_n d\nu\right|^p \right)^{1/p} +\left(\sum_{k=1}^N  \nu(A_k)^{-p/q} \left| \int_{U_k\setminus A_k} f_n d\nu\right|^p\right)^{1/p} \,\right]\\
        & \leq\ \liminf_{n\to\infty} \|f_n\|^p_{L^p(\nu)}\ +\ \liminf_{n\to\infty} \left( \sum_{k=1}^N \nu(A_k)^{-p/q} \cdot \|f_n\|_{L^p(\nu)}^p \,\nu(U_k\setminus A_k)^{p/q} \right)^{1/p}\\
        &<\ M\ +\ \eps.
    \end{align*}
    Since $\eps>0$ was taken arbitrarily and $f_*\geq 0$ $\nu_*$-a.e., we arrive at
    \[
        \sum_{k=1}^N \nu(A_k)^{1-p} \left|\int_{A_k} f_*d\nu_*\right|^p \ \leq\  \sum_{k=1}^N \nu(A_k)^{1-p} \left|\int_{U_k} f_*d\nu_*\right|^p \ \leq\ M^p.
    \]
    Thus, $f|_{K_m}\in L^p(S,\nu_*)$ for any $m\in\NN$ and $\|f_*|_{K_m}\|_p\leq M$. We conclude that $f_*\in L^p(S,\nu_*)$ and $\|f\|_p\leq M$. Finally, it is a simple observation that $f_*\ind_{S_*}\in L^p(S,\nu)$ and that $\|f_*\ind_{S_*}\|_{L^p(\nu)} = \|f_*\|_{L^p(\nu_*)}$.
\end{proof}

If the sequence of functions $f_n\in L^p_+(\nu)$ in Theorem \ref{Lem:BoundedWeakLimits} is uniformly bounded in $L^p$-norm, then $M=\liminf_{n\to\infty}\|f_n\|_p\leq\sup_{n}\|f_n\|_p<\infty$. Thus, we are motivated to define for any $M>0$ and $1\leq p\leq\infty$:
\[
    L^p_{+,M}(\nu) := \{f\in L^p_+(\nu)\colon \|f\|_p\leq M\}.
\]
We shall identify functions on $S$ that are in $L^1\cap L^p$ relative to a fixed positive Borel measure $\nu$ on $S$ with finite signed Borel measures on $S$ through the linear map $j_\nu:f\mapsto fd\nu$.  Additionally,  one has
\begin{clry}\label{clry:Lp-part closed in weak measure topo}
    Let $\nu$ be a positive Radon measure on $S$. Then, for any $1< p\leq \infty$ and $M>0$, $L^p_{+,M}(\nu)$ is closed in $\CM(S)$ for the weak topology on measures. Hence, it is complete in $\CM(S)_\FM$.
\end{clry}
\begin{proof}
    First, since the weak topology on $\CM^+(S)$ is metrizable (e.g. by the metric derived from the Fortet-Mourier norm) and $L^p_{+,M}(\nu)$, viewed as set of measures, consists of positive measures, the weak topology on measures is metrizable on $L^p_{+,M}(\nu)$. Thus, it suffices to show that for any sequence $(f_n)$ in $L^p_{+,M}(\nu)$ that weakly converges to a measure $\mu_*\in \CM^+(S)$, that $\mu_*\in L^p_{+,M}(\nu)$. This is established by Theorem \ref{Lem:BoundedWeakLimits}. The last statement follows from the first part, that $\|\cdot\|_\FM$ metrizes the weak topology on measures on $\CM^+(S)$ and the completeness of the latter for the $\|\cdot\|_\FM$-norm.
\end{proof}

\begin{remark}
    1.) The Radon-Nikodym Theorem is valid for $\sigma$-finite measures \cite[Theorem 4.2.2]{cohn2013measure}. If the reference measure is not $\sigma$-finite then the measure that is absolutely continuous with respect to the latter, may not have a density with respect to this measure. An example of such a situation is given e.g. in \cite[Exercise 3.2.13, p.92]{Folland:1999}. 
    Therefore, from the absolute continuity of $\mu_*$ with respect to the general Radon measure $\nu$ one cannot conclude by this theorem, that $\mu_*$ has a density with respect to $\nu$. The proof of this result in Theorem \ref{Lem:BoundedWeakLimits} is based on the absolute continuity of $\mu_*$ with respect to the restriction $\nu_*$ of $\nu$ to the $\sigma$-compact set $S_*$, which is $\sigma$-finite.\\
    2.) We do not wish to assume $\sigma$-finiteness of $\nu$ in Theorem \ref{Lem:BoundedWeakLimits}, because we want in the future to allow for $\nu$ being e.g. a Hausdorff measure or a Haar measure on a topological group, which may not be $\sigma$-finite. Lebesgue measure $\lambda$ on $\R^d$ is, of course.\\
    3.) The reader should notice that the lattice operations $\mu\mapsto \mu^+$ and $\mu\mapsto \mu^-$ are not continuous for the weak topology on measures. The results above therefore do not readily extend to sequences $(f_n)$ in $L^1\cap L^p(S,\nu)$ that do not consist of positive functions.\\
    4.) It is worth recalling the observation that Corollary \ref{clry:Lp-part closed in weak measure topo} cannot hold for $p=1$.
\end{remark}

\section{Invariance properties of the solution operator}
\label{invariance}

In this section, we begin proving the main invariance results.
We start with recalling or establishing essential estimates on the flow maps $\Phi^v_{t_0,t}$ in $\Omega$ and the associated Markov operators $P^v_{t_0,t}$ on $\CM(\Omega)$.

\subsection{Basic properties of the flow semigroup}

First, in the `classical setting' that we have for the time-dependent velocity field $v_\bullet$, namely that (V1)--(V3) are satisfied (i.e., $v_\bullet\in C(\R_+,\BL(\Omega,\R^d))$), the associated flow maps $\Phi^v_{t_0,t}$ to ODE \eqref{Eq:FlowODE} are bijective, Lipschitz with estimate \eqref{eq:estimate Lipschitz const flow} for the Lipschitz constant, and having Lipschitz inverse.  Moreover, for any $t_0\geq 0$, $\tau>0$, $t\in I_\tau=[t_0, t_0+\tau]$ and $x\in\Omega$,
\begin{equation}\label{eq:estimate flow from start}
    d(\Phi^v_{t_0,t}(x),x) = \bigl| \Phi^v_{t_0,t}(x) - x\bigr| \leq \int_{t_0}^t \bigl |v_s(\Phi^v_{t_0,s}(x))\bigr|\,ds \leq C_{t_0,\tau}|t-t_0|\ \to\ 0\qquad\mbox{as}\ t\downarrow t_0.
\end{equation}
Here, one can take $C_{t_0,\tau}=\sup_{s\in I_\tau}\|v_s\|_\infty$.

\begin{lemma} \label{lem:strong continuity Pt}
    For any $t_0\geq 0$ and $\mu\in\CM(\Omega)$, consider the motion $\varphi:t\mapsto P_{t_0,t}^v\mu:[t_0,\infty)\to\CM(\Omega)$. The following statements hold:
    \begin{enumerate}
        \item[(i)] $\varphi$ is continuous for the $C_b(\Omega)$-weak topology on measures.
        \item[(ii)] $\varphi$ is locally Lipschitz continuous for the $\|\cdot\|_\FM^*$-norm topology. In particular, for any bounded interval $I\subset [t_0,\infty)$, there exists a constant $C'_I>0$ such that 
        \begin{equation}\label{eq:local Lipschitz estimate MO flow}
            \|P^v_{t_0,t}\mu - P^v_{t_0,s}\mu\|_\FM^*\ \leq\ C'_I \|\mu\|_\TV\cdot |t-s|\qquad \mbox{for all}\ s,t\in I,\ \mu\in\CM(\Omega),
        \end{equation}
        where the constant $C'_I$ may be taken as
        \begin{equation}\label{eq:constant estimate MO flow}
            C'_I = \sup_{t\in I}\|v_t\|_\infty \cdot \exp\bigl(\int_{t_0}^{\bar{t}} |v_\tau|_L\,d\tau\bigr) = L^v_{t_0,\bar{t}}\, \sup_{t\in I}\|v_t\|_\infty,\qquad \bar{t}:=\sup(I).
        \end{equation}
    \end{enumerate} 
\end{lemma}
\begin{proof}
    {\it(i).}\ For any $t\geq t_0$ and $f\in C_b(\Omega)$,
    \[
        \langle P^v_{t_0,t}\mu, f\rangle = \langle\mu, f\circ \Phi_{t_0,t}^v\rangle = \int_\Omega f(\Phi^v_{t_0,t}(x))\,d\mu(x).
    \]
    For fixed $x\in \Omega$, $t\mapsto f(\Phi_{t_0,t}^v(x))$ is continuous. Lebesgue's Dominated Convergence Theorem gives the continuity of $t\mapsto \langle P^v_{t_0,t}\mu, f\rangle$.\\
    {\it(ii).}\ Take $f\in\BL(S)$ with $\|f\|_\infty\leq 1$ and $|f|_L\leq 1$. Let $I$ be a bounded interval in $[t_0,\infty)$ and $s,t\in I$. Then, using \eqref{eq:estimate Lipschitz const flow} and \eqref{eq:estimate flow from start}, one gets for any $\mu\in\CM(\Omega)$:
    \begin{align*}
        \bigl |\langle P^v_{t_0,t}\mu-P^v_{t_0,s}\mu, f\rangle \bigr|\  &\leq\ \int_\Omega \bigl| f(\Phi_{t_0,t}^v(x)) - f(\Phi^v_{t_0,s}(x))\bigr|\,d|\mu|(x)\\
        &\leq\ |f|_L |\Phi_{t_0,(s\wedge t)}|_L \int_\Omega d\bigl( \Phi^v_{(s\wedge t), (s\vee t)}(x),x\bigr) \,d|\mu|(x)\\
        &\leq\ |f|_L \cdot \exp \left( \int_{t_0}^{(s\wedge t)} |v_\tau|_L\,d\tau\right)\cdot C_I \bigl|(s\vee t)-(s\wedge t)\bigr| \cdot \|\mu\|_\TV.
    \end{align*}
    Here, one can take $C_I = \sup_{s\in I}\|v_s\|_\infty$, according to \eqref{eq:estimate flow from start}.
    Taking the supremum over $f$ in the unit ball in $\BL(\Omega)$ yields \eqref{eq:local Lipschitz estimate MO flow}, 
    where one can take $C'_I$ as in \eqref{eq:constant estimate MO flow}.
\end{proof}

We shall further need the following transport formula for measures that are absolutely continuous with respect to Borel-Lebesgue measure $\lambda$ on $\Omega$.
\begin{lemma}\label{lem:estimates push-forward flow}
    Let $u_0\in L^1(\Omega,\lambda)$ and $v_\bullet\in C(\R_+, \BL(\Omega,\R^d))$ be a bounded Lipschitz time-dependent vector field on $\Omega$. Then, for $t\geq t_0\geq 0$,
    \begin{equation}\label{eq:expression ut}
        P^v_{t_0,t}(u_0 \cdot\lambda) = u_t\cdot\lambda,\quad \mbox{with}\quad
        u_t(x) = u_0((\Phi^v_{t_0,t})^{-1}(x))\,\bigl| \det D\Phi^v_{t_0,t}\bigl((\Phi^v_{t_0,t})^{-1}(x)\bigr)\bigr|^{-1}.
    \end{equation}
    Moreover, for any $t\geq t_0\geq 0$ for $\lambda$-a.e. $x\in\Omega$,
    \begin{equation}\label{eq:estimate Jacobian flow}
        \exp\bigl(-\int_{t_0}^t\bigl\|[\mathrm{div} (v_s)]^-\bigr\|_\infty\,ds\bigr)\ \leq\ \det D\Phi^v_{t_0,t}(x)\ \leq\ \exp\bigl(-\int_{t_0}^t \bigl\|[\mathrm{div} (v_s)]^+\bigr\|_\infty\,ds\bigr).
    \end{equation}
\end{lemma}
\begin{proof} 
    The proof is based on the classical area formula for Lipschitz maps, see \cite{Ambrosio-Crippa:2014} p.1196-1197 for details, equation (2.5) in particular. Estimate \eqref{eq:estimate Jacobian flow} is inequality (2.6) in \cite{Ambrosio-Crippa:2014}, p.1197.
\end{proof}

\begin{remark}
    In various works one has studied flows defined by vector field that are less regular than Lipschitz, motivated -- among others -- by Fluid Mechanics and Control Theory, see e.g. \cite{DiPerna-Lions:1989,Bogachev-Mayer-Wolf:1999,Ambrosio-Crippa:2014}. Typically, Lebesgue measure is taken as a reference measure, like we do here. However, Bogachev and Mayer-Wolf \cite{Bogachev-Mayer-Wolf:1999} provide conditions under which the push-forward of a finite measure under the flow results into finite measures that are still absolutely continuous with respect to that measure. Although their results apply to probability measures as reference, instead of the Borel-Lebesgue measure, results similar to those stated above were obtained (cf. in particular \cite{Bogachev-Mayer-Wolf:1999}, Theorem 2.6, p.13). 
\end{remark}

\noindent The explicit expression for the push-forward under the flow and estimates in Lemma \ref{lem:estimates push-forward flow} gives 

\begin{prop}\label{clry:invaraiance Lp under Pt}
    For any $1\leq p\leq\infty$, and any $t\geq t_0\geq 0$,  $P^v_{t_0,t}$ leaves $L^p_+(\lambda)$ invariant. If $u_0\in L^p_+(\lambda)$, then $P^v_{t_0,t}(u_0\,d\lambda):=u_t\,d\lambda$, with $u_t$ given by \eqref{eq:expression ut}. If $q$ is the exponent conjugate to $p$, $\frac{1}{p}+\frac{1}{q}=1$, then  
    \begin{equation}\label{eq:estimate Lp norm ut}
        \|u_t\|_{L^p(\lambda)} \leq \|u_0\|_{L^p(\lambda)}\,\exp\left(\mbox{$\frac{1}{q}$}\int_{t_0}^t  \bigl\|[\mathrm{div}(v_s)]^-\bigr\|_\infty\,ds\right) \qquad\mbox{for all}\ t\geq 0.
    \end{equation}
\end{prop}

\begin{proof}
    Using that the flow operators $\Phi^v_{t_0,t}$ are Lipschitz, invertible, with Lipschitz inverse $(\Phi^v_{t_0,t})^{-1}$ and -- again -- the transformation formula for the Borel-Lebesgue measure, we arrive at
    \begin{align*}
        \int_\Omega u_t(x)^p\,dx \ & =\ \int_\Omega u_0\bigl( (\Phi^v_{t_0,t})^{-1}(x)\bigr)^p \bigl|\det D\Phi^v_{t_0,t}\bigl( (\Phi^v_{t_0,t})^{-1}(x) \bigr)\bigr|^{-p}\,dx\\
        &= \int_\Omega u_0(y)^p\bigl|\det D\Phi^v_{t_0,t}(y)\bigr|^{-p}\cdot 
        \bigl|\det D(\Phi^v_{t_0,t})^{-1}\bigl( \Phi^v_{t_0,t}(y) \bigr)\bigr|^{-1}\,dy
    \end{align*}
    The identity $(\Phi^v_{t_0,t})^{-1}\circ\Phi_{t_0,t}^v = \mathrm{Id}$ yields $\bigl|\det D(\Phi^v_{t_0,t})^{-1}(\Phi^v_{t_0,t}(y))\bigr|^{-1} = \bigl| \det D\Phi^v_{t_0,t}(y) \bigr|$ for $\lambda$-a.e. $y$. Application of \eqref{eq:estimate Jacobian flow} gives the result.
\end{proof}

\subsection{Equivalent formulation of mild solution by dilation}

Central in the argumentation is the following observation, inspired by \cite{Sikic} -- there employed for semigroups instead of flows -- and which appears also in the particular setting of multiple of the identity (but in $L^1$ instead of $\CM(\Omega)$) in \cite{Canizo_ea:2012}, proof of Lemma 3.2.

\begin{prop}\label{prop:equivalent equations mild solution}
    Let $\mu_\bullet :[t_0,T]\to\CM(\Omega)_\FM$ be a continuous and $\|\cdot\|_\TV$-bounded map and let $\nu\in\CM(\Omega)$. Let $B$ be a bounded linear operator on $\CM(\Omega)_\FM$ that commutes with all $P^v_{s,t}$, $s,t\in[t_0,T],\ s\leq t$.
    Then the following statements are equivalent:
    \begin{enumerate}
        \item[(i)] $\mu_\bullet$ is a mild solution to Problem (P) on $[t_0,T]$ with initial condition $\nu$.
        \item[(ii)] $\mu_\bullet$ satisfies, as equation in $\CM(\Omega)_\FM$:
        \begin{equation}\label{eq:VoC with perturb B}
            \mu_t\ = \ e^{-B(t-t_0)}P^v_{t_0,t}\nu\ +\ \int_{t_0}^t e^{-B(t-s)} P^v_{s,t}\bigl[ f_s(\mu_s) + B\mu_s \bigr]\,ds\quad \mbox{for all}\ t\in [t_0,T].
        \end{equation}
    \end{enumerate}  
    In particular, any mild solution $\mu_\bullet$ to Problem (P) on $[t_0,T]$ with $\mu_{t_0}=\nu$ satisfies for any $c\in\R$:
    \begin{equation}\label{eq:semilinear with c}
        \mu_t = e^{-c(t-t_0)}P^v_{t_0,t}\nu\ + \ \int_{t_0}^t e^{-c(t-s)}P^v_{s,t}\bigl[ f_s(\mu_s)+c\mu_s\bigr]\,ds\qquad\quad \mbox{for all}\ t\in[t_0,T].
    \end{equation}
\end{prop}
\begin{proof}
    The last statement immediately follows from {\it (ii)}, taking $B=cI$. So, consider the equivalences. First observe the following two simplifications of repeated integrals, which hold for any $k\in\NN$:
    \begin{equation}
        \int_{t_0}^t\int_{t_0}^s h(\sigma)(t-s)^{k-1}\,d\sigma\, ds\ =\ \int_{t_0}^t \int_{\sigma}^t h(\sigma)(t-s)^{k-1}\,ds\,d\sigma\ =\ \int_{t_0}^t h(\sigma)\,\mbox{$\frac{1}{k}$} (t-\sigma)^k\, d\sigma,
        \label{eq:simplifying a double integral}
    \end{equation}
    for a measurable function $h$ for which Fubini's Theorem holds.\\
    {\it (i)\ $\Rightarrow$\ (ii)}.\ Using that the operator $B$ commutes with all $P^v_{s,t}$, one can rewrite \eqref{eq:VoC formula mild solution} into
    \begin{equation}\label{eq:VoC mild sol shifted by B}
        \mu_t\ =\ P^v_{t_0,t}\nu\ +\ \int_{t_0}^t P^v_{s,t}\bigl[ f_s(\mu_s)+B\mu_s \bigr]\,ds\ -\ B\int_{t_0}^t P^v_{s,t}\mu_s\,ds.
    \end{equation}
    Substituting the right hand side of expression \eqref{eq:VoC mild sol shifted by B} for $\mu_s$ in the last integral in this equality and using \eqref{eq:simplifying a double integral} to simplify double integrals, one obtains, by repeating this procedure $n$ times:
    \begin{align*}
        \mu_t\ &=\ \sum_{k=0}^n \left[ (-1)^k\mbox{$\frac{(t-t_0)^k}{k!}$}\, B^k P^v_{t_0,t}\nu\ +\ \int_{t_0}^t  (-1)^k\mbox{$\frac{(t-s)^k} {k!}$}  B^kP^v_{s,t} \bigl[ f_s(\mu_s)+ B\mu_s\bigr]\,ds \right]\\
        &\qquad +\ (-1)^{n+1}B^{n+1}\int_{t_0}^t \mbox{$\frac{(t-s)^n} {n!}$}  P^v_{s,t} \mu_s\,ds.
    \end{align*}
    This equality can easily be verified by induction. Letting $n\to\infty$, one obtains {\it (ii)} in the limit in $\CM(\Omega)_\FM$, because $B$ is $\|\cdot\|_\FM^*$-bounded and $\mu_\bullet$ is $\|\cdot\|_\TV$-bounded.\\
    {\it (ii)\ $\Rightarrow$\ (i)}.\ Let $\mu_\bullet$ satisfy \eqref{eq:VoC with perturb B}. We must show that \eqref{eq:VoC formula mild solution} holds. We proceed -- again -- by repeatedly substituting equation \eqref{eq:VoC with perturb B} into itself and the equations thus obtained, and using \eqref{eq:simplifying a double integral} to simplify. Thus, after $n$ successive substitutions, we claim that one arrives at the equality
    \begin{align}
        \mu_t\ &=\ \sum_{k=0}^n \left[ \mbox{$\frac{(t-t_0)^k}{k!}$} B^k e^{-B(t-t_0)}P^v_{t_0,t}\nu\ +\ B^k \int_{t_0}^t \mbox{$\frac{(t-s)^k}{k!}$} e^{-B(t-s)} P^v_{s,t}\bigl[f_s(\mu_s)\bigr]\, \,ds \right] \nonumber\\
        &\qquad\ +\ B^{n+1} \int_{t_0}^t \mbox{$\frac{(t-s)^n}{n!}$} e^{-B(t-s)} P^v_{s,t} \mu_s \,\, ds. \label{eq:modified VoC n iterations}
    \end{align}
    Clearly, the equality holds for $n=0$, by assumption. Then, substituting \eqref{eq:VoC with perturb B} into \eqref{eq:modified VoC n iterations} and using that $B$ commutes with the operators $P^v_{s,t}$, one gets
    \begin{align*}
        \int_{t_0}^t e^{-B(t-s)} &P^v_{s,t} \mu_s \,\mbox{$\frac{(t-s)^n}{n!}$}\, ds\ =\ \int_{t_0}^t \mbox{$\frac{(t-s)^n}{n!}$} e^{-B(t-s)} P^v_{s,t}\, e^{-B(s-t_0)}P^v_{t_0,s}\nu\, \, ds \\
        &+\ \int_{t_0}^t \mbox{$\frac{(t-s)^n}{n!}$} e^{-B(t-s)} P^v_{s,t}\,\left[ \int_{t_0}^s e^{-B(s-\sigma)} P^v_{\sigma,s}\bigl[ f_\sigma(\mu_\sigma) \bigr]\, d\sigma\right]\, ds\\
        &\qquad +\ B \int_{t_0}^t \mbox{$\frac{(t-s)^n}{n!}$} e^{-B(t-s)} P^v_{s,t} \left[ \int_{t_0}^s e^{-B(s-\sigma)} P^v_{\sigma,s}\mu_\sigma \, d\sigma \right] \,ds\\
        & =\ \mbox{$\frac{(t-t_0)^{n+1}}{(n+1)!}$}\, e^{-B(t-t_0)} P^v_{t_0,s}\nu \ +\ \int_{t_0}^t \int_{t_0}^s \mbox{$\frac{(t-s)^n}{n!}$} e^{-B(t-\sigma)} P^v_{\sigma, t} \bigl[ f_\sigma(\mu_\sigma) \bigr]\, d\sigma\, ds\\
        &\qquad +\ B \int_{t_0}^t \int_{t_0}^s \mbox{$\frac{(t-s)^n}{n!}$} e^{-B(t-\sigma)} P^v_{\sigma, t}\mu_\sigma \,d\sigma\,ds\\
        & = \ \mbox{$\frac{(t-t_0)^{n+1}}{(n+1)!}$}\, e^{-B(t-t_0)} P^v_{t_0,s}\nu \ +\ \int_{t_0}^t \mbox{$\frac{(t-\sigma)^{n+1}}{(n+1)!}$} e^{-B(t-\sigma)} P^v_{\sigma, t} \bigl[ f_\sigma(\mu_\sigma) \bigr]\, d\sigma \\
        &\qquad +\ B \int_{t_0}^t \mbox{$\frac{(t-\sigma)^{n+1}}{(n+1)!}$} e^{-B(t-\sigma)} P^v_{\sigma, t}\mu_\sigma \,d\sigma.
    \end{align*}
    Thus, the induction step is completed and we conclude that \eqref{eq:modified VoC n iterations} holds for all $n\in\NN$. By letting $n\to\infty$ in \eqref{eq:modified VoC n iterations}, we obtain \eqref{eq:VoC formula mild solution}.  Mild solutions are unique. Therefore, $\mu_\bullet$ must be the mild solution to Problem (P) on $[t_0,T]$ with initial condition $\nu$.
\end{proof}

The equivalence of Problem (P) to the `dilated' problem \eqref{eq:VoC with perturb B} or \eqref{eq:semilinear with c} as observed in Proposition \ref{prop:equivalent equations mild solution} has a counterpart with regard to the Picard iteration operators $\Psi_{t_0,\tau}$ and those connected to the $B$-dilated problem \eqref{eq:VoC with perturb B}:
\begin{equation}
    \Psi^B_{t_0,\tau}(\nu_\bullet)(t)\ :=\ e^{-B(t-t_0)}P^v_{t_0,t}\nu_{t_0}\ +\ \int_{t_0}^t e^{-B(t-s)} P^v_{s,t}\bigl[ f_s(\nu_s) + B\nu_s\bigr]\,ds\quad \mbox{for}\ t\in [t_0,t_0+\tau].
\end{equation}
It is the following result:
\begin{lemma}\label{lem:contractive Picard B-dilated}
    Let $0\leq t_0<T$ and let $B$ be a bounded linear operator on $\CM(\Omega)_\FM$ that commutes with all $P^v_{s,t}$ for $s,t\in[t_0,T]$, $s\leq t$. Let $\delta>0$ and $\nu_0\in\CM(\Omega)_\FM$. Suppose that $t_1\in [t_0,T)$ and $\tau_1>0$ are such that $t_1+\tau_1\leq T$ and $\Psi_{t_1,\tau_1}$ is a strict contraction on $C_{\nu_0,\delta}([t_1,t_1+\tau_1],\CM(\Omega)_\FM)$. Let $\mu_\bullet^*$ be the unique fixed point in this space. Then, for any $r'>\|\nu_0\|_\TV+\delta$ there exists $N\in\NN$ and an equidistant partition $t_1=t'_0<t'_1<\dots<t'_N=t_1+\tau_1$ of $[t_1,t_1+\tau_1]$ of width $\tau':=\tau_1/N$, such that $\Psi^B_{t'_j,\tau'}$ is a strict contraction in $C_{\nu_j,\delta_j}([t'_j,t'_{j+1}],\CM(\Omega)_\FM)$ for all $j=0,1,\dots, N-1$, where $\nu_j := \mu^*_{t'_j}$ and $\delta_j := r' - \|\nu_j\|_\TV>0$. Moreover, $\nu_{j+1}$ is the end point at $t=t'_{j+1}$ of the unique fixed point of $\Psi^B_{t'_j,\tau'}$.
\end{lemma}

\begin{proof}
    Let $b$ be the operator norm of $B$ acting in $\CM(\Omega)_\FM$. If $b=0$, there is nothing to prove. So, assume $b>0$. Put $r:=\|\nu_0\|_\TV+\delta$. Note that $\|\mu^*_t\|_\TV\leq r$ for all $t\in[t_1,t_1+\tau_1]$. Let $r'>r$. For any $t_2\in[t_1,t_1+\tau_1)$ and $\tau_2>0$ such that $t_2+\tau_2\leq t_1+\tau_1$ and for any $\delta'>0$ and $\nu'\in \CM(\Omega)_\FM$ such that $\|\nu'\|_\TV+\delta'\leq r'$, one gets for $\nu^1_\bullet, \nu^2_\bullet\in C_{\nu',\delta'}([t_2,t_2+\tau_2],\CM(\Omega)_\FM)$ and $t\in[t_2,t_2+\tau_2]$ the estimates
    \begin{align*}
        \bigl\| \Psi^B_{t_2,\tau_2}(\nu^1_\bullet)(t) - \Psi^B_{t_2,\tau_2}(\nu^2_\bullet)(t) \bigr\|_\FM^*\ &\leq\ \int_{t_2}^t e^{-b(t-s)} \bigl\| P^v_{s,t}\bigl[ f_s(\nu^1_s) - f_s(\nu^2_s)\ +\ B(\nu^1_s-\nu^2_s)\bigr] \bigr\|_\FM^*\,ds\\
        & \leq \ \int_{t_2}^t e^{-b(t-s)} L^v_{s,t}( L^f_{r'} + b)\,\|\nu^1_s-\nu^2_s\|_\FM^*\, ds\\
        & \ \leq \ (L^f_{r'} + b) L^v_{t_1,t_1+\tau_1}\,\frac{1-e^{-b\tau_2}}{b}\, \|\nu^1_\bullet - \nu^2_\bullet\|_\infty.
    \end{align*}
    Pick $N\in\NN$ such that 
    \begin{equation}\label{eq:cond N}
        (L^f_{r'} + b) L^v_{t_1,t_1+\tau_1}\,\frac{1-e^{-b\tau_1/N}}{b}\ <\ 1.
    \end{equation}
    If $N=1$, the operator $\Psi^B_{t_1,\tau_1}$ is a strict contraction and the proof is complete. If $N\geq 2$, then put $\tau':=\tau_1/N$ and let $t'_j:= t_1+ j\tau'$ for $j=0,1,\dots, N-1$. Define $\nu_j:=\mu_{t'_j}^*$ and $\delta_j:= r' - \|\nu_j\|_\TV$. Since $\|\mu_t^*\|_\TV\leq r<r'$ for all $t\in[t_1,t_1+\tau_1]$, $\delta_j>0$ for all $j$. The assumptions and choices made yield the conclusion that $\Psi^B_{t'_j,\tau'}$ is a strict contraction on $C_{\nu_j,\delta_j}([t'_j,t'_{j+1}],\CM(\Omega)_\FM)$.\\
    The unique fixed point of $\Psi^B_{t'_j,\tau'}$, say $\mu^{(j)}_\bullet$, must be a mild solution to Problem (P) on $[t'_j,t'_{j+1}]$, according to Proposition \ref{prop:equivalent equations mild solution}. Mild solutions are unique (see Lemma \ref{lem:uniqueness mild solutions}). Thus, $\mu^{(j)}_\bullet$ equals $\mu^*$ on $[t'_j,t'_{j+1}]$. This yields that $\nu_{j+1}=\mu^{(j)}_{t'_{j+1}}$.   
\end{proof}

\subsection{Preservation of positivity}

Positivity of the mild solution to Problem (P) for positive initial conditions has been considered in \cite{ackleh2020well} Corollary 7.2 under the assumption that the reaction term $f$ can be written as a rate, or in \cite{Canizo_ea:2012} Lemma 3.2 for a particular type of equation. In the latter, the problem is reduced to positivity of the solutions to the PDE in $L^1(\Omega,\lambda)$, by a density argument. Below, however, we shall give an argument for positivity of solutions that is intrinsic to the measure setting, in the sense that it exploits the Banach lattice structure of $\CM(\Omega)$. Moreover, it allows to consider reaction terms that are more general than those considered in \cite{ackleh2020well}.

Proposition \ref{prop:equivalent equations mild solution} will be crucial in establishing the positivity of the mild solutions for positive initial conditions.
Assume that $f_\bullet$ satisfies (A1)--(A3) and let $\mu_\bullet$ be a mild solution to Problem (P) on $[t_0,T]$, for initial condition $\nu$. Then, $s\mapsto f_s(\mu_s)$ is continuous. In particular, it is strongly Bochner measurable. According to Lemma \ref{lem:measurability postive part}, for any $c\in\R$ the functions $s\mapsto \bigl[ f_s(\mu_s)+c\mu_s\bigr]^\pm$ are strongly measurable. Thus, the Bochner integral in $\overline{\CM}(S)_\BL$ in \eqref{eq:semilinear with c} can be split into one over the positive part and one over the negative part of $f_s(\mu_s)+c\mu_s$. Because of positivity of the operators $e^{-c(t-s)}P^v_{s,t}$ and closedness of the positive cone $\CM^+(\Omega)$ in $\overline{\CM}(\Omega)_\FM$, each integral separately takes value in $\CM^+(\Omega)$, while $e^{-c(t-t_0)}P^v_{t_0,t}\nu\in\CM^+(\Omega)$ as well for $\nu\in\CM^+(\Omega)$. We conclude from Proposition \ref{prop:equivalent equations mild solution} that for any $c\in \R$ the following inequality holds: 
\begin{equation}\label{eq:general estimate negative part}
    \mu^-_t \ \leq\ \int_{t_0}^t e^{-c(t-s)}P^v_{s,t}\bigl[ f_s(\mu_s) + c\mu_s\bigr]^-\,ds\quad\mbox{for all}\ t\in[t_0,T].
\end{equation}

This general observation allows drawing conclusions on the positivity of mild solutions for positive initial conditions for quite a general class of reaction terms $f_\bullet$.  The next result is a straightforward consequence of inequality \eqref{eq:general estimate negative part}:
\begin{prop}
    Let $f_\bullet$ satisfy (A1)--(A3) and let $t_0>0$. Let $\nu\in\CM^+(\Omega)$ and let $\mu_\bullet$ be the mild solution to problem (P) with initial condition $\mu_{t_0}=\nu$. Let $[t_0,T^*)$ be its maximal domain of existence. Assume that for any $R>0$ and $T>0$ there exists $c=c_{R,T}\in\R$ such that
    \begin{equation}\label{eq:conditin shifted positivity}
        f_t(\mu) + c\mu\ \geq\ 0,\qquad \mbox{for all}\ \mu\in\CM_R(\Omega),\ t\in[0,T].
    \end{equation}
    Then, $\mu_t\in\CM^+(\Omega)$ for all $t\in [t_0,T^*)$.
\end{prop}

\begin{proof}
    Let $T\in(t_0,T^*)$ and put $R:=\sup_{t\in[t_0,T]}\|\mu_t\|_\TV$. For $c= c_{R,T}$ as indicated, $[f_s(\mu_s)+c\mu_s]^-=0$ for all $s\in [t_0,T]$. So, \eqref{eq:general estimate negative part} implies that $\mu^-_t=0$ for $t\in [t_0,T]$.
\end{proof}

However, condition \eqref{eq:conditin shifted positivity} is too restrictive and often impractical for many applications. Interestingly, it can be improved especially if the term $f_\bullet$ has a specific shape. In such case, sufficient conditions ensuring the positivity of the solution can be reached so that they fit better to applications. The second result addresses this issue:

\begin{prop}\label{prop:positivity mild solution semilinear}
    Let $f_\bullet$ be of the form \eqref{eq:specific admiss non-linearity} in Proposition \ref{prop:example family admissible reaction terms}, such that it satisfies (A1)--(A3).
    Assume that $F_\bullet$ is such that for every $R>0$ and $T>0$, there exists $C_{R,T}>0$ such that  
    \begin{equation}\label{cond:positivity c}
        \sup_{0\leq t\leq T}\|F_t(\mu)\|_\infty\leq C_{R,T} \qquad \mbox{for all}\ \mu\in\CM_R(\Omega).
    \end{equation}
    Then, the mild solution $\mu_\bullet$ to Problem (P) with initial condition $\nu\in\CM^+(\Omega)$ is positive for all its time of existence.
\end{prop}
\begin{proof}
    Let $t_0\geq 0$ and $\nu\in\CM^+(\Omega)$. Let $T^*>t_0$ be the maximal time of existence of the mild solution $\mu_\bullet$ of Problem (P) with initial condition $\mu_{t_0}=\nu$. Take $T\in(t_0,T^*)$ and let $R:=\sup_{t\in[t_0,T]}\|\mu_t\|_\TV$.
    For any $c\in\R$ and $\mu\in\CM(\Omega)$ one has
    \begin{align}
        \bigl[f_t(\mu)+c\mu\big]^-\ &=\ \bigl[p_t(\mu) + \bigl(F_t(\mu)+c\bigr)\mu\bigr]^- \nonumber\\
        &\leq \ \bigl(F_t(\mu) + c\bigr)^- \mu^+ \ + \bigl(F_t(\mu) + c\bigr)^+\mu^-.\label{eq:estimate neg part}
    \end{align}
    Then, for $c\geq C_{R,T}$
    \[
        \bigl[f_t(\mu_t)+c\mu_t\big]^-\ \leq (C_{R,T} + c) \mu^-_t\qquad \mbox{for all}\ t\in[t_0,T]
    \]
    and, according to \eqref{eq:general estimate negative part}, we have
    \[
        0\ \leq \mu^-_t \leq\ (C_{R,T}+c)\int_{t_0}^t  e^{-c(t-s)} P^v_{s,t}\mu_s^-\,ds\qquad \mbox{for all}\ t\in[t_0,T].
    \]
    Note that the Bochner measurability of the integrand follows from Lemma \ref{lem:measurability postive part}. Taking the $\|\cdot\|_\TV$-norm on both sides, the ordering is preserved and the direct application of Proposition \ref{prop:measure Bochner integral and TV-norm} yields
    \[
        \|\mu^-_t\|_\TV\ \leq\ (C_{R,T}+c) \int_{t_0}^t e^{-c(t-s)} \|\mu^-_t\|_\TV\,ds\qquad \mbox{for all}\ t\in[t_0,T].
    \]
    Using Gr\"onwall's Inequality, one concludes that $\mu^-_t=0$ for all $t\in[t_0,T]$. 
\end{proof}

\begin{remark}
A similar conclusion as in Proposition \ref{prop:positivity mild solution semilinear} is reached  in \cite{ackleh2020well} Corollary 7.2 , but a (positive) production term $p_t(\mu)$ in the non-linearity $f_t(\mu)$ is not allowed by the assumptions made there. Moreover, the argument follows different lines of reasoning than ours.
\end{remark}

The result just proven uses the `global equation' \eqref{eq:VoC formula mild solution} that the mild solution satisfies. Realizing that the mild solution is locally obtained by Picard iteration, one may weaken condition \eqref{cond:positivity c} on $F_\bullet$ to the extent that the bound needs to apply only to {\it positive} measures. Thus, the `behaviour' of $F_\bullet$ on non-positive measures is irrelevant to the preservation of positivity. This weakening of condition goes to a cost: the need of a more subtle argumentation in the proof.

\begin{theorem}\label{thrm:positivity}
    Let $f_\bullet$ be of the form \eqref{eq:specific admiss non-linearity} in Proposition \ref{prop:example family admissible reaction terms}.
    Assume that $F_\bullet$ is such that for every $R>0$ and $T>0$, there exists $C'_{R,T}>0$ such that  
    \begin{equation}\label{cond:positivity condition 2}
        \sup_{0\leq t\leq T}\|F_t(\mu)\|_\infty\leq C'_{R,T} \qquad \mbox{for all}\ \mu\in\CM^+(\Omega)\ \mbox{with}\ \|\mu\|_\TV\leq R.
    \end{equation}
    Then, the mild solution $\mu_\bullet$ to Problem (P) with initial condition $\mu_{t_0}=\nu\in\CM^+(\Omega)$ is positive for all its time of existence.
\end{theorem}

\begin{proof}
    Let $[t_0,T^*)$ be the time domain of the maximal mild solution. Recall that $\mu_\bullet$ is obtained by piecing together the unique fixed points of the `Picard iteration operators' $\Psi_{t_i,\tau_i}$ (see \eqref{eq:Picard iteration operator}) for various times $t_i\in[t_0,T^*)$ and durations $\tau_i>0$. The fixed points are obtained as limit of the iterations of the operators $\Psi_{t_i,\tau_i}$, starting from the function $t\mapsto P^v_{t_i,t}\nu_i$ on $[t_i,t_i+\tau_i]$ (see also Remark \ref{rem:solution locally Picard iteration}). Each of the iterations and the limit will have total variation norm uniformly bounded by $R_i>0$, say, on this interval. 
    
    In view of Proposition \ref{prop:equivalent equations mild solution} and Lemma \ref{lem:contractive Picard B-dilated}, the mild solution can be locally obtained equally well by iteration of suitable $B$-dilated operators $\Psi^{B_i}_{t'_j,\tau'_i}$, e.g. with $B_i=c_iI$, $c_i\in\R$ and starting from a function $t\mapsto e^{-c_i(t-t'_j)}P^v_{t'_j,t}\nu'_j$, for suitable $\nu'_j$, on an equally spaced partition $t_i = t'_0<t'_1<\dots<t'_N=t_i+\tau_i$. The number $N$ of subintervals depends on $|c_i|$ and on $R_i$ in so far that it can be arranged that all iterations under $\Psi^{c_i I}_{t'_j,\tau'_i}$ -- and the limit -- have total variation uniformly bounded by $2R_i$ on each of the subintervals $[t'_j,t'_{j+1}]$ (see Lemma \ref{lem:contractive Picard B-dilated} and its proof. In particular, $N$ is conditioned by equation \eqref{eq:cond N}, with $b=|c_i|$ and $r'=2R_i$). 
    
    Fix $c_i\geq C'_{2R_i,t_i+\tau_i}$ and determine $N_i$ and a partition of $[t_i,t_i+\tau_i]$ as above, such that each $\Psi^{c_i I}_{t'_j,\tau'_i}$, $j=1,\dots, N_i$, is a strict contraction. Assume that $\nu'_j$ is positive. Then $\nu^{(0)}_t:=e^{-c_i(t-t'_j)}P^v_{t'_j,\tau'_i}\nu'_j$ is positive for every $t\in[t'_j,t'_j+\tau'_i]$. Define iteratively 
    \[
        \nu_\bullet^{(n+1)} := \Psi^{c_i I}_{t'_j,\tau'_i}(\nu^{(n)}_\bullet).
    \]
    One has, by construction, that $\|\nu_t^{(n)}\|_\TV\leq 2R_i$ for all $n$.
    Hence, if $\nu^{(n)}_\bullet$ is positive, then condition \eqref{cond:positivity condition 2} and application of \eqref{eq:estimate neg part} gives
    \[
        \bigl[ f_t\bigl(\nu^{(n)}_t\bigr) + c\nu^{(n)}_t \bigr]^-\ =\ 0.
    \]
    So, we conclude that $\nu^{(n)}_\bullet$ must be positive on $[t'_j,t'_j+\tau'_i]$. Since $\CM^+(\Omega)$ is closed in $\CM(\Omega)_\FM$, we find that the limit of the iterations $\nu^{(n)}_\bullet$ must be positive. 
    Hence, if $\mu_{t_i}$ is positive, then $\mu_\bullet$ is positive on $[t_i,t_i+\tau_i]$. Since $\nu\in\CM^+(\Omega)$, we conclude that $\mu_\bullet$ is positive on its maximal domain of existence.
\end{proof}

\subsection{Preservation of $L^p(\lambda)$-regularity}

We shall now consider an additional condition on the reaction term $f_\bullet$ that ensures that initial conditions that have $L^p$ density with respect to Borel-Lebesgue measure $\lambda$ on $\Omega$ for some $1< p\leq\infty$, lead to mild solutions to Problem (P) that are absolutely continuous with respect to $\lambda$  with the density functions remaining in $L^p(\lambda)$ for all time up to time of possible blow-up in $L^p$-norm of the density. It should be noted, that this time may generally be earlier than the maximal existence time of the mild solution in the space of measures. 

To avoid cumbersome notation, we write $\|\mu\|_{L^p(\lambda)}$ to mean the $L^p$-norm of the density of $\mu$ relative to $\lambda$. By convention, it is $\infty$, if $\mu$ is not absolutely continuous with respect to $\lambda$, or if the density is not an $L^p$-function.
\vskip 2mm

Let $1<p\leq\infty$. For our purpose, it is necessary to assume that for all $t\geq 0$ the reaction $f_t$ maps $L^1\cap L^p(\lambda)$ into itself, viewed as subspaces of $\CM(\Omega)$. With slight abuse of wording, we shall say in that case that `$f_t$ {\it maps $L^p(\lambda)$ into itself}'. To have our method of proof work, we shall additionally require an assumption similar to (A4) in $L^p$-setting. That is, $f_t$ must be an {\it $L^p$-bounded map, locally uniformly in time}. More precisely, 
\begin{quote}
    \begin{enumerate}
        \item[(A4$_p$)] {\it For any $t>t'$ and $r>0$, there exists $B^f_{t',t}(r)\geq 0$ such that
        \[
            \bigl\| f_s(\phi\, d\lambda) \bigr\|_{L^p(\lambda)}\ \leq\ B^f_{t',t}(r)\quad\mbox{for all}\ \phi\in L^p(\lambda),\ \|\phi\|_{L^p(\lambda)} \leq r,\ t'\leq s\leq t.
        \]}
    \end{enumerate}
    \end{quote}
must hold. Put
\[
    D^{v,p}_{s,t} := \exp\left( \mbox{$\frac{1}{q}$}\int_s^t \bigl\| [\mathrm{div}_x(v_\tau)]^-\bigr\|_\infty\,d\tau \right),\qquad 0\leq s\leq t,\ \mbox{$\frac{1}{p}+\frac{1}{q}=1$}.
\]
Our major result on `$L^p$-invariance' is then:
\begin{theorem}\label{thrm:Lp invariance positive}
    Let $f_\bullet$ be of the form \eqref{eq:specific admiss non-linearity} and satisfy the conditions of Theorem \ref{thrm:positivity}. Moreover, assume that for some $1<p\leq\infty$, $f_t$ maps $L^p(\lambda)$ into itself for all $t\geq 0$ and is an $L^p$-bounded map, locally uniformly in time, i.e. satisfies (A4$_p$). Let $\nu=\phi d\lambda\in \CM^+(\Omega)$ be an initial condition at time $t_0$ with $\phi\in L^p_+(\lambda)$. Then the following statements holds:
    \begin{enumerate}
        \item[(i)] There exists $\hat{T}^*_\nu\leq T^*_\nu$ such that the maximal mild solution $\mu_\bullet$ to Problem (P) in $\CM(\Omega)$ takes values $\mu_t\in L^p_+(\lambda)$ for all $t\in[t_0,\hat{T}^*_\nu)$.
        \item[(ii)] For every $T\in (t_0,\hat{T}^*_\nu)$, $\sup_{t_0\leq t\leq T} \|\mu_t\|_{L^p(\lambda)}<\infty$.
        \item[(iii)] If $T^*_\nu<T^*_\nu$, then $\limsup_{t\to \hat{T}^*_\nu} \|\mu_t\|_{L^p(\lambda)}=+\infty$.
    \end{enumerate} 
\end{theorem}

\noindent Before getting to the details of the proof, we like to make several observations.

\begin{remark}
    1.)\ Theorem \ref{thrm:Lp invariance positive} requires the $L^p$-initial condition to be positive and requires the assumption of the conditions of Theorem \ref{thrm:positivity} to assure that the mild solution to Problem (P) remains positive throughout its maximal domain of existence. Positivity is required for application of the limit results, Theorem \ref{Lem:BoundedWeakLimits} and corollary, that allow to conclude that `$L^p$-ness' is preserved in the weak limit of measures.\\
    2.) In many applications, the measure-valued solution will represent a distribution of individuals or `mass' in space. It is then natural to have a model that preserves positivity of the solution. Moreover, the shape \eqref{eq:specific admiss non-linearity} for the reaction map $f_\bullet$ is natural as well, from modelling perspective. It has at time $t$ a production term $p_t$ and a degradation term that is absolutely continuous with respect to the current distribution, $F_t(\mu_t)\cdot\mu_t$, which is in accordance with the principle that `mass can only disappear where there is mass present'.\\
    3.) In specific cases, it may happen in principle that $L^p$-regularity of the measure-valued solution is lost at time $\hat{T}^*_\nu$ strictly before time $T^*_\nu$, while the former continues to exist as measure with less regular density or as a measure that may even be singular with respect to Lebesgue measure.  See e.g. \cite{Schmeiser, Velazquez, Winkler}, although these concern a different type of equation. 
\end{remark}

\begin{proof}
    {\it (Theorem \ref{thrm:Lp invariance positive}).}\ Let $\mu_\bullet$ be the maximal mild solution to Problem (P) on $[t_0,T^*_\nu)$. According to Theorem \ref{thrm:positivity}, $\mu_t\in\CM^+(\Omega)$ for all $t$. Let $t'\in [t_0,T^*_\nu)$ such that $\nu':=\mu_{t'}$ has an $L^p$-density with respect to $\lambda$. Following the proof Theorem \ref{thrm:positivity}, there exist $c\geq 0$, $\delta'>0$ and $\tau'>0$ such that $t'+\tau'< T^*_\nu$ and the Picard iteration operator $\Psi^c_{t',\tau'}$ is a strict contraction on $C_{\nu',\delta'}([t',t'+\tau'], \CM(\Omega)_\FM)$. Moreover, on $[t', t'+\tau']$, $\mu_t = \lim_{n\to\infty} \mu^{(n)}_t$ in $\CM(\Omega)_\FM$, as in \eqref{eq:Picard iteration in measures}, with $\Psi_{t',\tau'}$ replaced by $\Psi^c_{t',\tau'}$, and $c$ is such that each $\mu^{(n)}_\bullet$ is positive.

    Put $D:=D^{v,p}_{t',t'+\tau'}$. One has for $t'\leq t\leq t'+\tau'$, using Proposition \ref{clry:invaraiance Lp under Pt}:
    \[
        \bigl\| \mu^{(0)}_t\bigr\|_{L^p(\lambda)}\ =\ 
        e^{-c(t-t')} \bigl\| P^v_{t',t}\nu' \bigr\|_{L^p(\lambda)}\ \leq\ D^{v,p}_{t',t} \|\nu'\|_{L^p(\lambda)}\ \leq D\|\nu'\|_{L^p(\lambda)}
    \]
    and
    \[
        \bigl\| \mu^{(n+1)}_t\bigr\|_{L^p(\lambda)}\ =\ 
        \bigl\| \Psi^c_{t',\tau'}\bigl( \mu^{(n)}_\bullet \bigr) \bigr\|_{L^p(\lambda)}
        \leq\bigl\| \mu^{(0)}_t\bigr\|_{L^p(\lambda)}+\left\| \int_{t'}^t e^{-c(t-s)} P^v_{s,t} \bigl[ f_s(\mu^{(n)}_s) + c\mu^{(n)}_s \bigr]\,ds\right\|_{L^p(\lambda)}.
    \]
    Let $I_n(t)$ denote the measure-valued Bochner integral in $\overline{\CM}(\Omega)_\FM$ in the last term. We shall show by induction that 
    \begin{equation}\label{eq:property choice hat-tau}
        \mbox{{\it There exists}}\ 0<\hat{\tau}\leq\tau',\ \mbox{{\it such that}}\ \bigl\| \mu^{(n)}_t\bigr\|_{L^p(\lambda)}\leq 2D\|\nu'\|_{L^p(\lambda)}\ \mbox{{\it for all}}\ n\in\NN_0,\ t\in[t',t'+\hat{\tau}].
    \end{equation}
    For $n=0$ it is clearly satisfied.  Take $r:=2D\|\nu'\|_{L^p(\lambda)}$ and put $B^f_r:= \sup_{t\in[t',t'+\tau']} B^f_{t',t}(r)$, with $B^f_{t',t}(r)$ as in Assumption (A4$_p$). According to \cite{Diestel-Uhl}, we have that
    \[
        I_n(t)\ \in\ (t-t') \cdot \overline{\mathrm{co}}\left( \bigl\{ e^{-c(t-s)} P^v_{s,t}\bigl[f_s(\mu^{(n)}_s) + c\mu^{(n)}_s \bigr]\colon t'\leq s\leq t\bigr\}  \right),
    \]
    where the closure of the convex hull is taken in $\overline{\CM}(\Omega)_\FM$. For $t'\leq s\leq t\leq t'+\tau'$, by induction hypothesis,
    \[
        \bigl\| e^{-c(t-s)} P^v_{s,t}\bigl[f_s(\mu^{(n)}_s) + c\mu^{(n)}_s \bigr] \bigl\|_{L^p(\lambda)}\ \leq\ e^{-c(t-s)} D^{v,p}_{s,t}(B_r^f + cr)\ \leq\ D(B^f_r + cr).
    \]
    Any point in the convex hull is positive by assumption on $c$ and will have the same bound for its $L^p$-norm. $I_n(t)$ is the limit in the weak topology on measures of $(t-t')$-multiple of a sequence of elements in the convex hull. Because of positivity,  Theorem \ref{Lem:BoundedWeakLimits} yields that $I_n(t)\in L^p_+(\lambda)$ and
    \begin{equation}
        \bigl\| I_n(t) \bigr \|_{L^p(\lambda)}\ \leq\ (t-t') D(B^f_r + cr)\quad\mbox{for all}\ n\in\NN_0.
    \end{equation}
    One can choose $0<\hat{\tau}\leq\tau'$, independent of $n$, such that
    \begin{equation}
        \hat{\tau}\cdot D(B_r^f + cr) \leq D\|\nu'\|_{L^p(\lambda)}.
    \end{equation}
    Then, \eqref{eq:property choice hat-tau} is satisfied.
    The iterates $\mu^{(n)}_\bullet$ being in $L^p_+(\lambda)$, uniformly bounded in $L^p$-norm, we can again invoke Theorem \ref{Lem:BoundedWeakLimits}, yielding $\mu_t\in L^p_+(\lambda)$ for all $t\in [t',t'+\hat{\tau}]$, with $\|\mu_t\|_{L^p(\lambda)}\leq 2D\|\nu'\|_{L^p(\lambda)}$.

    Now, the `classical extension argument' can be repeated for extending local solutions to a maximal one, starting from $\nu\in L^p_+(\lambda)$ at $t=t_0$. This then yields the wanted statements.
\end{proof}

\section{Concluding considerations and observations}

We now want to address some natural questions in relation to the invariance results for the solution operator that we have established in the previous section.
\vskip 2mm

We have shown in Theorem \ref{thrm:Lp invariance positive} that under fairly mild and common assumptions (compare e.g. with those in \cite{ackleh2020well}) a positive initial condition with $L^p$-density with respect to Lebesgue measure leads to a mild solution to Problem (P) that has a positive $L^p$-density for some positive amount of time. That may not be the full maximal existence time of the mild solution in the space of measures.

Mild solutions are by definition required to be continuous in time.
Thus, this concept allows to define a dynamical system on (suitable subsets) of measures, which have -- usually -- continuous motions. What can one say at this point about regularity in time of the mild solution with values in $L^p(\lambda)$? 
This remains elusive to us at this moment. For finite reference measure instead of the $\sigma$-finite Lebesgue measure, we have found ways to say more about this question. However, the step to non-finite measures is in this context not as easy to make as may be anticipated. This is a topic for further investigation in future works.

Finally, we like to point the reader to a particular `feature' of assumptions that were imposed. The Lipschitz continuity of the reaction term $f_\bullet$ in measures and in Fortet-Mourier norm, formulated as Assumption (A2), is providing the convergence of the Picard iterations in the space of measures with Fortet-Mourier norm. The condition (A4$_p$) then yields solutions $\mu_\bullet\in L^r([t_0,T], L^p(\lambda))$.  To get an $L^p$-valued solution, one `usually' imposes on the reaction map an assumption of Lipschitz type in $L^p(\lambda)$. Here it is merely $L^p$-boundedness together with a Lipschitz condition in Fortet-Mourier norm. On the other hand, we cannot establish any regularity in time, continuity in particular, that is obtained from an $L^p$-Lipschitz condition.

Moreover, there is no clear relationship known (yet) between $\|\mu\|_{L^p(\lambda)}$ and $\|\mu\|_\FM^*$ for a signed measure $\mu\in L^1\cap L^p(\lambda)$, as occurs in computing distances. Thus, a Lipschitz condition on the reaction map in either norm cannot be related.

\begin{appendix}
\section{Bochner measurability of lattice operations on measures}

Let $S$ be a Polish space and $d$ an admissible metric on $S$. It should be noted that the maps $\mu\mapsto\mu^\pm$, where $\mu^+$ and $\mu^-$ are the positive and negative part of $\mu$ in the Hahn-Jordan decomposition, are not continuous for the weak topology on measures. However, one has the following result -- and a consequence. As far as we are aware, they could not be found in the literature and are of separate interest in (Bochner) integration of measure-valued functions.

\begin{lemma}\label{lem:Jordan is measurable}
    For every $f\in B_b(S)$, $\mu\mapsto \langle\mu^+, f\rangle = \int_S f\,d\mu^+$ is Borel measurable as map from $\CM(S)_\FM$ to $\R$. In particular, if $E$ is a Borel set in $S$, then $\mu\mapsto \mu^\pm(E)$ and $\mu\mapsto \mu(E)$ are Borel measurable maps $\CM(S)_\FM\to\R$.
\end{lemma}

\begin{proof}
The proof is similar to the one of Lemma 3.2.6, p. 15, in \cite{vanDenzen}; it uses the Monotone Class Theorem for functions.
\end{proof}

It implies the following result concerning strongly measurable functions in Bochner integration. A function with values in a Banach space is strongly measurable if it is the point-wise limit of a sequence of simple step functions. 
Let $(S',\Sigma')$ be a measurable space.
\begin{lemma}\label{lem:measurability postive part}
    Let $f:S'\to \CM(S)_\FM$ be strongly measurable. Then $f^\pm: S'\to \CM(S)_\FM:s'\mapsto [f(s')]^\pm$ are strongly measurable.
\end{lemma}
\begin{proof}
We use Pettis' Measurability Theorem (see, e.g.,  \cite{Diestel-Uhl}). Since $\CM(S)_\FM$ is separable, it suffices to show that for every $\phi\in\CM(S)_\FM^*$, $\phi\circ f^+$ is Borel measurable. Observe that $\CM(S)_\FM^*\simeq \BL(S)$: $\phi\in\CM(S)_\FM^*$ corresponds to $g_\phi\in\BL(S)$ through $\phi(\mu) = \int_S g\,d\mu$ (cf. \cite{Hille-Worm}). Since $g_\phi$ is measurable, there exist simple step functions
\[
g_n = \sum_{k=1}^{N_k} a^{(k)}_n 1_{E^{(k)}_n}
\]
such that $g_n\to g_\phi$ point-wise while $|g_n(x)|\leq |g_\phi(x)|$ for all $x\in S'$ (see \cite{cohn2013measure}). Then for any $s'\in S'$,
\begin{equation}\label{eq:comp as limit}
    \phi\circ f^+(s') = \lim_{n\to\infty} \int_S g_nd\bigl[ f^+(s')\bigr] = \lim_{n\to\infty} \sum_{k=1}^{N_n} a^{(k)}_n [f(s')]^+\bigl(E^{(k)}_n\bigr).
\end{equation}
Because $s'\mapsto f(s'):S'\to\CM(S)_\FM$ and $\mu\mapsto \mu^+(E^{(k)}_n)$ are Borel measurable (Lemma \ref{lem:Jordan is measurable}), every term in the sum in \eqref{eq:comp as limit} is a Borel measurable function. So the point-wise limit $\phi\circ f^+$ is too.
\end{proof}

\end{appendix}

\section*{Acknowledgments}
S. C. H. and A. M. thank G. S. Magnuson Foundation for a partial financial support via the grant nr. MG2020-0008.


\end{document}